\documentclass[10pt,amssymb]{article}
\usepackage{amsfonts,amssymb,amsmath,latexsym,verbatim,amscd}

\title{Foliations by constant mean curvature tubes}

\author{Rafe Mazzeo \\ Stanford University \and Frank Pacard \\
Universit\'e Paris XII}

\date{}

\newtheorem{theorem}{Theorem}[section]
\newtheorem{proposition}{Proposition}[section]
\newtheorem{corollary}{Corollary}[section]
\newtheorem{lemma}{Lemma}[section]
\newtheorem{definition}{Definition}[section]
\newtheorem{remark}{Remark}[section]

\newcommand{\R}{\mathbb{R}}
\newcommand{\N}{\mathbb{N}}

\newcommand{\Z}{\mathbb{Z}}
\newcommand{\e}{\varepsilon}
\newcommand{\del}{\partial}
\newcommand{\ds}{}

\newcommand{\calC}{{\mathcal C}}
\newcommand{\calH}{{\mathcal H}}
\newcommand{\calL}{{\mathcal L}}

\newcommand{\calO}{{\mathcal O}}
\newcommand{\calT}{{\mathcal T}}
\newcommand{\calU}{{\mathcal U}}

\newcommand{\frakJ}{{\mathfrak J}}

\begin{document}

\maketitle

\begin{abstract}
Let $\Gamma$ be a nondegenerate geodesic in a compact Riemannian
manifold $M$. We prove the existence of a partial foliation of
a neighbourhood of $\Gamma$ by CMC surfaces which are small
perturbations of the geodesic tubes about $\Gamma$. There
are gaps in this foliation, which correspond to a bifurcation 
phenomenon. Conversely, we also prove, under certain restrictions, 
that the existence of a partial CMC foliation of this type about 
a submanifold $\Gamma$ of any dimension implies that $\Gamma$ is minimal. 
\end{abstract}

\section{Introduction}
Constant mean curvature hypersurfaces constitute a very important
class of submanifolds in a compact Riemannian manifold
$(M^{n+1},g)$. In this paper we are interested in families of such
submanifolds, with mean curvature varying from one member of the
family to another, which form (partial) foliations and which
`condense' to a submanifold $\Gamma \subset M$ of codimension
greater than $1$.
Our main results concern the existence of such families and,
conversely, the geometric nature  of the submanifolds $\Gamma$ to
which such families can condense.

The simplest case, where $\Gamma$ is a point, was considered by Ye
a decade ago, \cite{Ye-1}, \cite{Ye-2}. He proved that if $p \in
M$ is a nondegenerate critical point of the scalar curvature
function $R_g$, then there exists a neighborhood $\calU \ni p$
such that $\calU \setminus \{p\}$ is foliated by constant mean
curvature (for short CMC) spheres; in fact, the members of this
family are small perturbations of the geodesic spheres of radius
$\rho$, $0 < \rho < \rho_0$, and hence they have mean curvatures
$H = 1/\rho \to \infty$. Moreover, this foliation is essentially
unique. Conversely, if a neighbourhood of $p$ admits such a
foliation, then necessarily $\left. \nabla R_g \right|_p = 0$. In
very closely related work, Ye \cite{Ye-3}, and by quite different
methods (using inverse mean curvature flow) Huisken and Yau
\cite{Hui-Yau}, proved the existence of a unique foliation by CMC
spheres near infinity in an asymptotically flat manifold (of
nonnegative scalar curvature); this is of interest in general
relativity.

In this paper we study the existence of families of CMC
hypersurfaces which converge to a (closed, embedded)  submanifold
$\Gamma^\ell \subset M^{n+1}$, particularly in the case $\ell=1$.
Define the geodesic tube
\[
\calT_\rho(\Gamma) := \{q \in M^{n+1}:
\quad \mbox{dist}_g (q,\Gamma)= \rho \};
\]
this is a smooth hypersurface provided $\rho$ is smaller than
the radius of curvature of $\Gamma$, and we henceforth always
tacitly assume that this is the case. The mean curvature of this tube satisfies
\begin{equation}
H_{\calT_\rho(\Gamma)} = \frac{n-\ell}{n \rho} + \calO(1) \qquad
\mbox{as}\qquad \rho \searrow 0, \label{eq:1-1}
\end{equation}
and hence it is plausible that we might be able to perturb this
tube to a CMC hypersurface with $H \equiv (n-\ell)/(n\rho)$. It
turns out that this may not be possible for every (small) $\rho >
0$,  or when $\ell > 1$, but we prove the~:
\begin{theorem} Suppose that $\Gamma$ is a simple closed embedded
geodesic with nondegenerate Jacobi operator. Then there exist $k_0
\in {\mathbb N}$ and sequences $\rho_k'< \rho_k'' \to 0$, for $k
\geq k_0$ such that when $\rho \in I_k := (\rho_k',\rho_k'')$, the
geodesic tube $\calT_\rho(\Gamma)$ may be perturbed to a CMC
hypersurface $\Sigma_\rho$ with $H = \frac{n-1}{n\rho}$. The
$\Sigma_\rho$ are nonintersecting and foliate the open set equal
to their union, hence they form a partial foliation of some
neighborhood of $\Gamma$. \label{th:existence}
\end{theorem}
The hypersurface $\Sigma_\rho$ is a small perturbation of
$\calT_\rho(\Gamma)$ in the sense that it is the normal graph (for
some function whose $L^\infty$ norm is bounded by a constant times
$\rho^3$) over a small translate of $\Gamma$ (by some translation
whose $L^\infty$ norm is bounded by a constant times $\rho^2$),
cf.\ \S 4 for the precise formulation. In addition, we have
rather precise information on the location and width of the intervals
$I_k$:
\begin{equation}
\begin{array}{rcl}
\rho_k' -  \frac{\sqrt{n-1}\,\Lambda}{2\,\pi\,(k+1)} & = &
 {\cal O} (k^{-9/4}),
\\[3mm]
 \rho_k''-\frac{\sqrt{n-1}\,\Lambda}{2\,\pi\,k}  & = &
 {\cal O} (k^{-9/4}),
\end{array}
\label{eq:rhok}
\end{equation}
where $\Lambda$ equals the length of $\Gamma$. Note that the
nondegeneracy condition on $\Gamma$ is a mild and generic one, and
that there are no stringent conditions on the curvature along
$\Gamma$, as is the case in Ye's theorem when $\Gamma$ is a point.

The existence of gaps in this foliation stems from the fact that
at certain radii, the Jacobi operators on the geodesic tubes
$\calT_\rho(\Gamma)$ become degenerate, and this substantially
complicates certain analytic steps in the construction. However,
this gap behaviour is a real phenomenon, and is linked to a bifurcation
phenomenon, as we explain more carefully in \S 5.

Recall that the index of a compact CMC or minimal submanifold is
the number of negative eigenvalues of (the negative of) its Jacobi
operator. We estimate the index of the leaves of the partial
foliation we construct~:
\begin{proposition}
If $\rho \in I_k$, then for $k$ sufficiently large,
\[
\mathrm{Index}\, (\Sigma_\rho)  =  \mathrm{Index}\, (\Gamma) + 2
\, k + 1.
\]
\label{pr:index}
\end{proposition}
In particular, $\mbox{Index}(\Sigma_\rho) \to \infty$ as $\rho \to
0$. It is easy to show that for a generic metric on $M$, the
moduli space of CMC hypersurfaces condensing to $\Gamma$ is a
smooth one-dimensional set, i.e. a (possibly infinite) union of
curves, and that the index is constant on each component. Thus
Proposition~\ref{pr:index} shows that the leaves $\Sigma_\rho$ for
$\rho \in I_k$ lie in different components of this moduli space.

The final part of this paper concerns necessary conditions on
$\Gamma$ in order that a sequence of CMC surfaces condensing to
$\Gamma$ exists. We show that if there exists a sequence of CMC
hypersurfaces $\Sigma_j$, each of which can be written (in an
appropriate sense) as a normal graph over the geodesic tube
$\calT_{\rho_j}(\Gamma)$ where $\rho_j \to 0$, then $\Gamma$ must
be minimal. We defer to \S 5.2 for the precise statement of this
result. We have proved this converse only under rather stringent
conditions, but posit the following:

\smallskip
\noindent{\bf Conjecture:}
{\it Let $\Gamma$ be a closed embedded $\ell$-dimensional submanifold
of $M$ and that there exist sequences $\rho_k' < \rho_k'' \to
0$ and a partial foliation by CMC hypersurfaces $\Sigma_\rho$,
$\rho \in I_k := (\rho_k',\rho_k'')$ (with $\rho_k'' < \rho_{k-1}'$),
satisfying:
\begin{itemize}
\item[(i)] The mean curvature of $\Sigma_\rho$ equals
$\frac{n-1}{n \rho}$;

\item[(ii)] There exists a constant $c >0$, independent of $k$ and
$\rho \in I_k$ such that
\[
\Sigma_\rho \subset  \{q \in M^{n+1}: \quad \mbox{dist}_g
(q,\Gamma) \leq c \, \rho \} ;
\]
\item[(iii)] The norm of the second fundamental form of these
hypersurfaces satisfy $|A_{\Sigma_\rho} |\leq c \,\frac{1}{\rho}$
for some constant $c >0$, again independent of $k$ and $\rho \in
I_k$.
\end{itemize}
Then $\Gamma$ is a minimal submanifold. }

\smallskip

One might even be able to weak hypotheses (ii) and (iii)
substantially, but even with these hypotheses, the proof is
already probably difficult. We have chosen to prove this converse
only under much stronger hypotheses in order to include one main
calculation which explains why the minimality of $\Gamma$ is the
natural conclusion.

On the other hand, our method of proof encounters serious analytic
difficulties when $\dim \Gamma > 1$, and it is unclear whether
there is a general result concerning existence of
families of CMC hypersurfaces concentrating along a minimal
submanifold of dimension greater than $1$. The technical
complications are due ultimately to the lack of sufficiently good
estimates for $(-\Delta_\Gamma - \lambda)^{-1}$ on H\"older spaces
when $\lambda \to \infty$ and $\lambda$ lies in a spectral gap.
(When $\Gamma$ is a curve, the spectral gaps are large and the
Jacobi operator is an ODE, and such estimates are easy to obtain!).

There are some parallels between our Theorem~1 and some recent
results concerning solutions of the equation
\[
\e^2 \, \Delta u + f(u) = 0
\]
with vanishing Neumann data on a smooth, bounded domain $\Omega
\subset \R^n$. For example, the second author and Ritor\'e
\cite{Pac-Rit} prove the existence of positive solutions to
\[
\e^2 \, \Delta u - u^3 + u = 0
\]
which concentrate along a minimal submanifold $\Sigma$ as $\e
\searrow 0$. On the other hand, Malchiodi and Montenegro
\cite{Mal-Mon} construct positive solutions of
\[
\e^2 \, \Delta u + u^3 - u = 0
\]
which concentrate along $\del \Omega$ as $\e \searrow 0$. As in
the present paper, the same `spectral gap' phenomenon limits their
results to domains $\Omega \subset {\mathbb R}^2$.

\section{Geometry of tubes}
In this section we derive expansions as $\rho \searrow 0$ for the
metric, second fundamental form and mean curvature of the tubes
$\calT_\rho(\Gamma)$ and their perturbations. There is a famous
and beautiful formula for the volume of these tubes, due
originally to Herman Weyl, when $M$ has constant curvature, which
has found applications in fields as diverse as geometric measure
theory and statistics. We refer to Gray's monograph \cite{Gra} for
Weyl's formula and references  to later work.

\subsection{Fermi coordinates and Taylor expansion of the metric
near $\Gamma$} We first consider the asymptotic development of the
metric $g$ in Fermi coordinates around $\Gamma$. This leads to an
asymptotic formula for the metric on the geodesic tubes
$\calT_\rho(\Gamma)$. These computations are standard, and are
described more systematically  in \cite{Gra}.

Fix an arclength parametrization $\gamma(t)$ of $\Gamma$, $t \in
[a,b] :=I$, and denote by  $SN\Gamma$ the sphere bundle in
$N\Gamma$. Then
\[
SN\Gamma \ni (t,v) \longmapsto \exp_{\gamma(t)}(\rho v) \in
\calT_\rho(\Gamma)
\]
is a diffeomorphism when $\rho$ is small enough. Choose a parallel
orthonormal frame $E_1, \cdots, E_n$ for $N\Gamma$ (along $(a,b)$,
say). This determines a coordinate system
\[
x := (x_0, x_1 \ldots, x_n) \longmapsto
\exp_{\gamma(x_0)}(x_1E_1 + \ldots + x_n E_n) := F(x),
\]
and the corresponding coordinate vector fields $X_\alpha : =  F_*
(\del_{x_\alpha})$. We write $x' = (x_1,\ldots,x_n)$, and adopt
the convention that indices $i,j,k, \ldots \in \{1,\ldots,
n\}$, whereas $\alpha, \beta, \ldots \in \{0,\ldots, n\}$.

\begin{remark}
For simplicity, we identify the metric $g$ on $M$ and its
pullback $F^*g$ on some neighbourhood in $\R \times \R^n$. This
allows us to use the linear operations on the latter space. With
slight abuse of notations, we identify $F(x_0, x')$ with $(x_0,
x')$ and $X_\alpha$ with $\del_{x_\alpha}$.
\end{remark}

We also use cylindrical Fermi coordinates. Thus let $r =
\sqrt{x_1^2 + \ldots + x_n^2}$, which by Gauss' lemma is the
geodesic distance from $x$ to $\Gamma$. The vector
\begin{equation}
\del_r =  \ds \frac{1}{r}\sum_{i=1}^{n} x_i \, X_i
\label{eq:defnu}
\end{equation}
is the unit normal to the geodesic tubes.

We have arranged that the metric coefficients $g_{\alpha \beta} =
\langle X_\alpha , X_\beta \rangle$ equal $\delta_{\alpha\beta}$
along $\Gamma$. We now compute higher terms in the Taylor
expansions of these functions. In the following, the notation
$\calO(r^m)$ indicates a function $f$ such that it and its partial
derivatives of any order, with respect to the vector fields $X_0$
and $x_i \, X_j$, are bounded by $C \, r^m$ in some fixed
$T_{\rho_0}(\Gamma)$. Also, we shall compute the metric
coefficients at a point $q: = F(x_0, x')$ in terms of geometric
data at $p : = F (x_0, 0)$ and the radius $r = d(p,q)$.

We begin with the expansion of the covariant derivative:
\begin{lemma}
For $\alpha,\beta = 0, \ldots, n$,
\begin{equation}
\nabla_{X_\alpha } \, X_\beta = \sum_{\gamma=0}^n {\cal O}(r)
X_\gamma, \label{eq:3-2a}
\end{equation}
and for $\alpha=\beta =0$, we record the more precise expansion
\begin{equation}
\nabla_{X_0} \, X_0 = -  \ds \sum_{i,j=1}^n \langle R( X_j, X_0)
\, X_i , X_0 \rangle_p  \, x_i \, X_j + \sum_{\gamma=0}^n {\cal O}
(r^2) \, X_\gamma. \label{eq:3-2b}
\end{equation}
\label{le:3-1}
\end{lemma}
\noindent{\bf Proof:}
Anywhere on $\Gamma$,
\[
\nabla_{X_0} X_0 = \nabla_{X_0} X_j = \nabla_{X_j} X_0 =
\nabla_{X_i} X_j = 0.
\]
The vanishing of the first two terms is obvious since $\Gamma$ is
a geodesic and the $X_i$ are parallel along it. Because we are
using coordinate vector fields, $\nabla_{X_\alpha}X_\beta =
\nabla_{X_\beta}X_\alpha$ for any $\alpha,\beta$, even away from
$\Gamma$, and this implies the vanishing of the third term. Since
any $X \in N_p\Gamma$ is tangent to the geodesic $\exp_p(sX)$, and
so $\nabla_{X_i+X_j}(X_i + X_j)=0$ at $p$, hence $\nabla_{X_i}X_j
+ \nabla_{X_j}X_i = 0$ there. Combined with the symmetry
statement, we obtain that the final term also vanishes. This now
gives (\ref{eq:3-2a}).

Next, using (\ref{eq:3-2a}) we get
\[
\begin{array}{rlllll}
X_i \langle \nabla_{X_0} X_0, X_j\rangle & = & \langle
\nabla_{X_i} \nabla_{X_0} X_0, X_j \rangle +
\langle \nabla_{X_0} X_0, \nabla_{X_i} X_j\rangle \\[3mm]
& = & \langle \nabla_{X_i} \nabla_{X_0} X_0, X_j \rangle +
\calO (r^2) \\[3mm]
& = & \langle R(X_i, X_0)\, X_0, X_j \rangle_p + \langle
\nabla_{X_0} \nabla_{X_i} X_0, X_j\rangle_p  + \calO(r^2)\\[3mm]
& = & \langle R(X_i, X_0)\, X_0, X_j\rangle_p + \calO(r)
\end{array}
\]
This implies (\ref{eq:3-2b}). \hfill $\Box$

\smallskip

Our next result gives the expansion of the metric coefficients in
Fermi coordinates. The expansion of the $g_{ij}$, $i,j=1, \ldots,
n$, agrees with the well known expansion for the metric in normal
coordinates, cf. \cite{Sch-Yau}, \cite{Lee-Par} or \cite{Will},
but we briefly recall the proof here for completeness.
\begin{proposition}
In the same notation as above, we have
\begin{equation}
\begin{array}{rllll}
g_{ij}(q) & = & \delta_{ij} + \frac{1}{3} \, \langle R (X_k
,X_i)\, X_\ell , X_j \rangle_p \, x_k \, x_\ell + {\cal O} (r^3)
\\[3mm]
g_{0i}(q) & = & {\cal O} (r^2) \\[3mm]
g_{00}(q) & = & 1 + \langle R (X_k , X_0 ) \, X_\ell , X_0
\rangle_p \,  x_k  \, x_\ell + {\cal O} (r^3).
\end{array}
\label{eq:3-1}
\end{equation}
\label{pr:3-1}
\end{proposition}
\noindent{\bf Proof:}
The function
\[
X_k \, g_{\alpha \beta} = \langle
\nabla_{X_k}X_\alpha,X_\beta\rangle + \langle X_\alpha,
\nabla_{X_k}X_\beta\rangle
\]
vanishes on $\Gamma$, and thus the first order terms vanish in all
of these Taylor expansions.

To compute the second order terms, it suffices to compute
\[
\begin{array}{rcl}
X_k \, X_k \, g_{\alpha\beta}(p) & = & X_k \, X_k \, \langle
X_\alpha,X_\beta \rangle
\\[3mm]
& = & \langle \nabla_{X_k}^2 X_\alpha,X_\beta\rangle +
\langle X_\alpha, \nabla_{X_k}^2 X_\beta\rangle +
2\langle \nabla_{X_k}X_\alpha,\nabla_{X_k}X_\beta\rangle
\end{array}
\]
and then polarize (i.e.\ replace $X_k$ by $X_k + X_\ell$, etc.).
By (\ref{eq:3-2a}), the final term vanishes. Also,
\[
\nabla_{X_k}^2 X_\alpha = \nabla_{X_k}\nabla_{X_\alpha}X_k =
\nabla_{X_\alpha} \nabla_{X_k}X_k + R(X_k,X_\alpha) \, X_k.
\]

First let $\alpha = j \geq 1$ and compute the first term on the
right. Since $\nabla_{X}\nabla_{X}X =0$ on $\Gamma$ for any $X$
which is a constant linear combination of the $X_i$, we have
\[
0 = \nabla_{X_k + \e X_j}\nabla_{X_k + \e X_j}(X_k + \e X_j);
\]
equating the coefficient of $\e$ to $0$ gives
$\nabla_{X_j}\nabla_{X_k}X_k = -2 \nabla_{X_k}\nabla_{X_k}X_j$,
and hence
\[
3 \, \nabla_{X_k}^2 X_j =  R(X_k,X_j) \, X_k,
\]
so finally
\[
X_k \, X_k \, g_{ij} = \frac{2}{3} \, \langle R(X_k, X_i) \, X_k.
X_j\rangle.
\]
The formula for the second order Taylor coefficient for $g_{ij}$
now follows at once.

When $\alpha=0$, $\nabla_{X_0} \nabla_{X_k}X_k \equiv 0$ on
$\Gamma$, so
\[
\nabla_{X_k}^2 X_0 =  R(X_k,X_0) \, X_k.
\]
from which it follows that
\[
X_k \, X_k \, g_{00} = 2 \, \langle R(X_k, X_0) \, X_k, X_0
\rangle
\]
and this gives the formula for $g_{00}$.

The second order Taylor coefficient for $g_{0i}$ has not been
given because it is not needed later. \hfill $\Box$

\subsection{Perturbed tubes and their mean curvature}

We now describe a suitable class of deformations of the geodesic
tubes $\calT_\rho(\Gamma)$, depending on a section $\Phi$ of
$N\Gamma$ and a scalar function $w$ on the spherical normal bundle
$SN\Gamma$. One of the main technical parts of this paper, which
occupies the rest of \S 2, is the computation of the mean
curvature of these hypersurfaces, at least asymptotically in
$\rho$ and for sufficiently small $\Phi$ and $w$.

The spherical normal bundle is locally trivialized by the map
\[
(a,b) \times S^{n-1} \ni (x_0,\theta) \longmapsto (\gamma(x_0),\sum \theta_j E_j)
\in SN\Gamma.
\]
Fix $\rho > 0$, and define
\[
G(x_0,\theta) :=  F\,\big(x_0,\rho \,(1+ w(x_0,\theta))\,\theta+\Phi(x_0)\big);
\]
the image of this map will be called $\calT_\rho(w,\Phi)$. Thus
$\calT_{\rho}(w,\Phi)$ is obtained by first taking the normal
graph of the function $\rho\, w$ over the tube of radius $\rho$ in
$N\Gamma$ and then translating by $\Phi$. In particular
\[
\calT_{\rho}(0,0) = \calT_\rho(\Gamma).
\]
It will sometimes be useful to calculate using a coordinate system
\[
\R^{n-1} \ni y \mapsto \Upsilon (y) \in S^{n-1},
\]
with associated coordinates vector fields $Y_j = \del_{y_j}\Upsilon$.
In particular, we regard $G$ as a function of $(x_0,y)$ and write
\[
G(x_0,y) :=  F \, \big(x_0,\rho \,(1+ w(x_0,y))\,\Upsilon+\Phi(x_0)\big).
\]

Two different types of H\"older spaces will be used to measure
regularity of functions on $SN\Gamma$ and sections of $N\Gamma$:
first, we use the ordinary H\"older spaces
$\calC^{m,\alpha}(SN\Gamma)$, $\calC^{m,\alpha}(\Gamma,N\Gamma)$,
but we shall also use modified H\"older spaces
$\calC^{m,\alpha}_\rho(SN\Gamma)$, $\calC^{m,\alpha}_\rho(\Gamma, N\Gamma)$
which are based on differentiations with respect to the vector fields
$\rho \, \del_{x_0}$ and $\del_{y_j}$ (where $y$ is any local coordinate
system on $S^{n-1}$, see above). Note that this is tantamount to using the
rescaled variable $s = x_0/\rho$ since $\del_s = \rho \,
\del_{x_0}$. We shall assume that
\[
\Phi(x_0) = \sum_{j=1}^n \phi_j(x_0) \, X_j \in \calC^{2,\alpha}(\Gamma,N\Gamma),
\qquad w \in \calC^{2,\alpha}_\rho(SN\Gamma).
\]

For $p \in \Gamma$, let $S_p$ denote the spherical fibre of $SN\Gamma$ over
$p$. Any function $w$ on $SN\Gamma$ decomposes into a sum of three terms
\[
w= w_0 + \hat{w} + \tilde{w},
\]
where the restriction to any $S_p$ of each of these terms lies in
the span of the eigenfunctions $\varphi_j(\theta)$ on $S^{n-1}$
with $j=0$, $j = 1, \ldots, n$, and $j > n$, respectively. The
first component, $w_0$, is a function on $\Gamma$ itself. Next,
the eigenfunctions $\varphi_j$, $1 \leq j \leq n$, are the restrictions
to $S^{n-1}$ of linear functions on $\R^n$, so any linear combination
of them can be identified with a translation in $\R^n$ (the linear function
$x \rightarrow a \cdot x$ being identified with the translation
$x \rightarrow x+a$). Correspondingly, the summand $\hat{w}$ is canonically
associated to a section $\Phi$ of the normal bundle $N\Gamma$.

We shall typically assume that the functions $w$ has `linear component' $\hat{w}
\equiv 0$, and shall regard the linear part of the perturbation as a section
of $N\Gamma$, as just described.

It will be fundamental in the analysis below to regard $w$ as a
function of $s = x_0/\rho$ and $y_j$, but $\Phi$ and $\gamma$ as
functions of $x_0$ (in particular, whenever we write  $\Phi'$, we
mean $\del_{x_0}\Phi := \sum \del_{x_0} \phi_j(x_0) \, X_j$).
However, we sometimes also write $G = G(s,y)$. For example, the
tangent space to $\calT_\rho(w, \Phi)$ is spanned by the vector
fields
\begin{equation}
\begin{array}{rcccl}
Z_0 & = & G_* (\del_s) & = &  \rho \, ( X_0 + \del_s w \, \Upsilon +
\Phi'),
\\[3mm]
Z_j &= &  G_* (\del_{y_j}) & = & \rho  \, ((1+w)\, Y_j +
\del_{y_j}w \, \Upsilon ), \qquad j=1, \ldots, n.
\end{array}
\label{eq:defz0zj}
\end{equation}

\begin{definition}
In the following, $L(w,\Phi)$ denotes any expression which is
a linear differential operator (of order at most $2$) in $w$ and $\Phi$
which satisfies
\begin{equation}
\|L(w, \Phi)\|_{{\mathcal C}_\rho^{0, \alpha}} \leq c \, \left( \|
w \|_{{\mathcal C}_\rho^{2, \alpha}(SN\Gamma) } +
\|\Phi\|_{{\mathcal C}^{2, \alpha}(\Gamma, N\Gamma) } \right),
\label{eq:AAA}
\end{equation}
where $c$ is independent of $\rho$. Similarly, $Q(w,\Phi)$ denotes any
nonlinear differential operator (of order less than or equal to
$2$) in $w$ and $\Phi$ which vanishes quadratically in the pair
$(w, \Phi)$ and such that
\begin{eqnarray}
\|Q(w_2, \Phi_2) - Q(w_1, \Phi_1) \|_{{\mathcal C}^{0,\alpha}_\rho}
& \leq & c \,\sup_{i=1,2} \left(\| w_i \|_{{\mathcal C}^{2,\alpha}_\rho (SN\Gamma)}
+ \|\Phi_i\|_{{\mathcal C}^{2,\alpha}(\Gamma, N\Gamma) } \right) \notag \\
\times \ \bigg(\| w_2-w_1 \|_{{\mathcal C}^{2,\alpha}_\rho
(SN\Gamma)} & + & \|\Phi_2  - \Phi_1 \|_{{\mathcal C}^{2,\alpha}
(\Gamma,N\Gamma) } \bigg) \label{eq:BBB}
\end{eqnarray}
Here the spaces ${\cal C}^{0, \alpha}_\rho$ are either equal to
${\cal C}^{0, \alpha}_\rho (SN\Gamma)$ or ${\cal C}^{0,
\alpha}_\rho (\Gamma, N \Gamma)$ according to the range of $L$ and
$Q$.  Finally, terms denoted $\calO(\rho^k)$ are bounded in
$\calC^{m,\alpha}(SN\Gamma)$ or $\calC^{m,\alpha}(\Gamma,
N\Gamma)$ by $C \,\rho^k$, where the constant $C$ does not depend
on $(w,\Phi)$ or $\rho$.
\end{definition}

\subsection{The first fundamental form}
The next step is the computation of the coefficients of the first
fundamental form of $\calT_\rho(w, \Phi)$ with respect to the
coordinates $(s,y)$. At the point
\[
q = F(\rho \, s ,\rho(1+w (s,y)) \Upsilon(y)+\Phi(\rho s))
\]
(and $p= F (\rho s,0)$), we obtain directly from (\ref{eq:3-1})
that
\begin{equation}
\begin{array}{rcl}
\langle X_0, X_0\rangle_q  & = & 1 + {\cal O} (\rho^2) + \rho \,
L(w,\Phi) + Q(w, \Phi) \\[3mm]
\langle X_i, X_j \rangle_q & = &  \langle X_i, X_j\rangle_p +
\frac{\rho^2}{3} \, \langle R(\Upsilon,X_i) \, \Upsilon,X_j\rangle_p + \calO
(\rho^3)  \\[3mm]
& {} & +\,  \frac{\rho}{3} \left[ \langle R(\Upsilon,X_i) \, \Phi ,X_j
\rangle_p + \langle R(\Phi,X_i) \, \Upsilon,X_j\rangle_p \right] \\[3mm]
& {} & +\, \rho^2\, L(w,\Phi) + Q(w,\Phi) \\[3mm]
\langle X_i, X_0 \rangle_q & = &  {\cal O} (\rho^2) + \rho \, L(w,
\Phi) + Q(w,\Phi).
\end{array}
\label{lems}
\end{equation}

We use these expansions to obtain the expansion of the first
fundamental form of ${\mathcal T}_\rho(\Phi,w)$.
\begin{proposition}
We have
\begin{equation}
\begin{array}{rlllll}
\ds \rho^{-2} \, \langle Z_0, Z_0 \rangle & = & 1 + \calO(\rho^2)
+ \rho \, L(w,\Phi) + Q(w,\Phi)\\[3mm]
\ds \rho^{-2} \, \langle Z_0, Z_j \rangle & = & \calO(\rho^2)
+ L(w,\Phi)+Q(w,\Phi) \\[3mm]
\ds \rho^{-2} \,  \langle Z_i, Z_j \rangle & = & \langle Y_i , Y_j
\rangle_p + \frac{\rho^2}{3}\,\langle R(\Upsilon,Y_i)\,
\Upsilon, Y_j \rangle_p  + \calO(\rho^3) \\[3mm]
& {} & +\, 2 \, w \, \langle Y_i,Y_j \rangle_p + \frac{\rho}{3} \,
\left[ \langle R(\Upsilon,Y_i) \, \Phi , Y_j \rangle_p   + \langle
R(\Upsilon,Y_j)\, \Phi , Y_i \rangle_p \right]  \\[3mm]
& {} & \, \rho^2 \, L(w,\Phi) + Q(w,\Phi).
\end{array}
\label{eq:3-4}
\end{equation}
\label{pr:3-2}
\end{proposition}
{\bf Proof :} The first equation is clear. We give more details about
how to derive the second and third estimates since the same argument will
be used frequently. First, it follows from (\ref{lems}) that
\[
\begin{array}{rcl}
\langle \Upsilon, Y_j \rangle_q & = &  \langle \Upsilon, Y_j \rangle_p +
\frac{\rho^2}{3} \, \langle R(\Upsilon,\Upsilon) \, \Upsilon,Y_j\rangle_p + \calO
(\rho^3) \\[3mm]
& {} & +\, \frac{\rho}{3} \left[ \langle R(\Upsilon,\Upsilon) \, \Phi ,Y_j
\rangle_p + \langle R(\Phi,\Upsilon) \, \Upsilon,Y_j\rangle_p \right] \\[3mm]
& {} & +\, \rho^2\, L(w,\Phi) + Q(w,\Phi)
\end{array}
\]
However, when $w=\Phi=0$, $\langle \Upsilon, Y_j \rangle_q =0$
since $\Upsilon$ is normal and $Y_j$ is tangent to ${\mathcal
T}_\rho (0,0)$ then, so that the sum of the first three terms on
the right, which is independent of $w$ and $\Phi$, must also vanish.
This, together with the fact that $R(\Upsilon, \Upsilon) =0$
implies that
\begin{equation}
\langle \Upsilon, Y_j \rangle_q = \frac{\rho}{3} \, \langle
R(\Phi,\Upsilon) \, \Upsilon,Y_j\rangle_p + \rho^2\, L(w,\Phi) +
Q(w,\Phi))
\label{eq:jjjjj}
\end{equation}
In particular, we get $\langle \Upsilon, Y_j \rangle_q  = \rho \,
L(w,\Phi) + Q(w,\Phi))$. The second and third equations follow
directly from this. \hfill $\Box$

\subsection{The normal vector field}
The next task is to find expansions for the unit normal to $\calT_\rho(w,\Phi)$.
We begin with the preparatory
\begin{lemma}
The following expansions hold
\[
\begin{array}{rllllll}
\langle \Upsilon, \Upsilon \rangle_q & = & 1 +  \rho^2 \, L(w, \Phi) + Q (w, \Phi)\\[3mm]
\langle \Upsilon , Z_0 \rangle_q & = & \rho \, L(w, \Phi) + \rho \, Q(w, \Phi)\\[3mm]
\langle \Upsilon, Z_j \rangle_q  & = & \rho \, \del_{y_j}w +
\frac{\rho^2}{3} \, \langle R( \Phi  , \Upsilon ) \, \Upsilon ,
Y_j \rangle_p + \rho^3 \, L(w, \Phi) + \rho \, Q (w, \Phi)
\end{array}
\]
\label{le:3-2}
\end{lemma}
{\bf Proof~:} These follow from (\ref{lems}). As at the end of the
last subsection, we are using that $\langle \Upsilon, \Upsilon
\rangle_q = 1$ when $w=\Phi=0$ and $R(\Upsilon , \Upsilon) =0$ to
obtain the first two expansion. The second expansion follows from
the fact that $\langle \Upsilon , Z_0 \rangle_q =0$ when
$w=\Phi=0$. Finally, to obtain the last expansion, we use
that $\langle \Upsilon, Z_j \rangle_q  = 0$ when $w=\Phi=0$
as well as the first expansion and (\ref{eq:jjjjj}). \hfill $\Box$

\smallskip

We can now proceed with the expansion of the unit normal vector
field to ${\mathcal T}_\rho (w, \Phi)$.
\begin{proposition}
The normal vector field $N$ to $\calT_\rho(w, \Phi)$ has the
expansion
\begin{equation}
\begin{array}{rllll}
N  & : = &  - \, \Upsilon + \sum_{j=1}^{n-1} \, \alpha_j
\, Y_j + \left( L(w,\Phi) + Q(w, \Phi) \right) \, X_0 \\[3mm]
& {} & +\, \sum_{j=1}^{n-1} \left( \rho^2 \, L(w, \Phi) + Q(w,\Phi)\right) \, X_j
\end{array}
\label{eq:3-5}
\end{equation}
where the coefficients $\alpha_j$ are solutions of the system
\[
\sum_{j=1}^{n-1} \alpha_j \, \langle Y_j, Y_i\rangle_p =
\del_{y_i} w  + \frac{\rho}{3} \langle R (\Phi, \Upsilon) \,
\Upsilon, Y_i \rangle_p  .
\]
\label{pr:3-3}
\end{proposition}
{\bf Proof~:} Define the vector field
\[
\tilde N  : =  - \, \Upsilon + a_0 \, Z_0 +  \sum_{j=1}^{n-1} \, a_j \,Z_j,
\]
by choosing the coefficients $a_\alpha$ so that that $\tilde N$ is
orthogonal to all of the $Z_\alpha$.

It follows at once from Lemma~\ref{le:3-2} and (\ref{eq:3-4})
that $\rho \, a_\alpha =  L(w, \Phi) + Q(w, \Phi)$ for every $\alpha$.
Plugging this back into each of the equations $\langle N, Z_i
\rangle_q = 0$ (thus neglecting the orthogonality condition when
$\alpha=0$ now), we find that the $a_j$ are solutions of the system
\[
\sum_{j=1}^{n-1} a_j \, \langle Y_j, Y_i\rangle_p =  \frac{1}{
\rho} \, \del_{y_i} w  + \frac{1}{3} \langle R (\Phi, \Upsilon) \, \Upsilon, Y_i
\rangle_p  + \rho  \, L(w, \Phi) + \frac{1}{ \rho} \, Q(w, \Phi).
\]
Recall also that $Z_j = \rho Y_j + \rho L(w,\Phi)$ so that $a_i Z_i =
\rho \alpha_i Y_i + \rho L(w,\Phi)$. Finally, we have
\[
|\tilde N|_q  =  1 + \rho^2 \, L(w, \Phi) + Q(w, \Phi).
\]
This gives (\ref{eq:3-5}). \hfill $\Box$

\subsection{The second fundamental form}
The most arduous step is the computation of the second fundamental
form. To simplify the computations below, we henceforth assume
that, at the point $\Upsilon(y) \in S^{n-1}$,
\begin{equation}
\langle Y_i, Y_j\rangle_p =\delta_{ij} \qquad \mbox{and}\qquad
\overline \nabla_{Y_i} Y_{j} =0 \label{eq:3-55}, \quad i,j = 1, \ldots, n-1
\end{equation}
(where $\overline \nabla$ is the connection on $TS^{n-1}$).
\begin{proposition}
The following expansions hold
\begin{equation}
\begin{array}{rllll}
\rho^{-2} \, \langle N, \nabla_{Z_0} Z_0 \rangle_q & = & \rho \,
\langle R(\Upsilon,X_0) \, \Upsilon , X_0 \rangle_p +
{\cal O} (\rho^2) \\[3mm]
& {} & - \,\frac{1}{\rho} \, \del_s^2 w  - \langle \Phi'', \Upsilon \rangle_p
+ \langle R ( \Upsilon , X_0) \, \Phi , X_0 \rangle_p  \\[3mm]
& {} & \, \rho \, L(w, \Phi) + \frac{1}{\rho} \, Q(w, \Phi),
\end{array}
\label{eq:3-6a}
\end{equation}
\begin{equation}
\rho^{-2} \, \langle N, \nabla_{Z_0 } Z_j \rangle_q = {\cal
O}(\rho) + \frac{1}{\rho} \, L (w, \Phi)  + \frac{1}{\rho}\,
Q(w,\Phi), \label{eq:3-6b}
\end{equation}
\begin{equation}
\begin{array}{rllll}
\rho^{-2} \, \langle N, \nabla_{Z_j} Z_{j} \rangle_q & = &
\frac{1}{\rho} + \frac{2}{3} \, \rho \, \langle R(\Upsilon, Y_j )
\, \Upsilon, Y_{j} \rangle_p +  {\cal O} ( \rho^2) \\[3mm]
& {} & -\, \frac{1}{\rho}  \del_{y_j}^2 w  + \frac{1}{\rho} \, w  +
\frac{2}{3} \, \langle R(\Phi, Y_j) \, \Upsilon , Y_j \rangle_p  \\[3mm]
& {} & +\, \rho  \, L (w, \Phi)  +  \frac{1}{\rho} \, Q (w, \Phi)
\end{array}
\label{eq:3-6c}
\end{equation}
\begin{equation}
\begin{array}{rllll}
\rho^{-2} \, \langle N, \nabla_{Z_i} Z_{j} \rangle_q & = & {\cal
O} (\rho) + \frac{1}{\rho}  \, L (w, \Phi)  + \frac{1}{\rho} \, Q
(w, \Phi), \quad i \neq j.
\end{array}\label{eq:3-6d}
\end{equation}
\label{pr:3-4}
\end{proposition}
{\bf Proof~:} First note that by Lemma~\ref{le:3-1}
\begin{equation}
\left. \nabla_{X_\alpha} \, X_\beta\right|_q = \sum_{\gamma=0}^n \left( {\cal O}
(\rho) + L(w, \Phi)+ Q(w, \Phi)\right)  \, X_\gamma,
\label{eq:3-7}
\end{equation}
since the coordinates of $q$ depend on $w$ and $\Phi$. Hence , as
\begin{equation}
\nabla_{Z_\alpha}  X_\beta  =  \sum_{\gamma=0}^n \left({\cal O}
(\rho^2) + \rho \, L(w, \Phi)+ \rho \, Q(w, \Phi)\right)  \,
X_\gamma, \label{eq:3-9}
\end{equation}
which follows from (\ref{eq:3-7}) and the fact that $Z_\alpha =
\rho \, \sum_\gamma (1+ L(w, \Phi) ) \, X_\gamma$.

We will also use that
\begin{equation}
N +  \Upsilon = \sum_{\alpha =0}^n \left(L(w, \Phi)+ Q(w,
\Phi)\right) \, X_\alpha, \label{eq:3-8}
\end{equation}
which follows from (\ref{eq:3-5}). Finally, we will need the
expansions
\begin{equation}
\begin{array}{rllllll}
\langle \Upsilon , X_0 \rangle_q & = & \rho \, L(w, \Phi) +  Q(w, \Phi)\\[3mm]
\langle \Upsilon, Y_j \rangle_q  & = & \rho \, L(w, \Phi) + Q (w,
\Phi)
\end{array}
\label{eq:3-10}
\end{equation}
whose proof can be obtained, as in Lemma~\ref{le:3-2}, starting
from (\ref{lems}).

\smallskip

\noindent
{\bf Estimate (\ref{eq:3-6a}):} We must expand
\[
\rho^{-2} \, \langle N, \nabla_{Z_0} Z_0 \rangle_q = \rho^{-1} \,
\left( \langle N, \nabla_{Z_0} X_0 \rangle_q + \langle N,
\nabla_{Z_0} (\del_s w \, \Upsilon) \rangle_q + \langle N,
\nabla_{Z_0} \Phi' \rangle_q \right)
\]
The estimate is broken into three steps:

\smallskip

\noindent {\bf Step 1} From (\ref{eq:3-5}) and Lemma~\ref{le:3-2} we get
\[
\begin{array}{rlllll}
\langle N, \Upsilon \rangle_q & = & - \langle \Upsilon, \Upsilon
\rangle_q + \sum_{j=1}^n \alpha_j \,\langle Y_j, \Upsilon\rangle_q
+  (L(w, \Phi) + Q(w,\Phi) ) \, \langle X_0 ,  \Upsilon \rangle_q \\[3mm]
& {} & +\, \sum_{j=1}^n \, (\rho^2 \, L(w, \Phi) + Q(w, \Phi) ) \,
\langle X_j, \Upsilon \rangle_q  \\[3mm]
& = & -  1 + \rho^2 \, L(w, \Phi) + Q(w, \Phi)
\end{array}
\]
Substituting $N = -\Upsilon + N+\Upsilon$ gives
\[
\langle N, \nabla_{Z_0} \Upsilon \rangle_q  = -  \frac{1}{2} \, \del_s
\langle \Upsilon, \Upsilon \rangle_q + \langle N + \Upsilon, \nabla_{Z_0} \Upsilon \rangle_q;
\]
by Lemma~\ref{le:3-2}
\[
\del_s \langle \Upsilon, \Upsilon \rangle_q = \rho^2 \, L(w, \Phi)+  Q(w, \Phi),
\]
and (\ref{eq:3-8}) and (\ref{eq:3-9}) imply
\[
\langle N + \Upsilon , \nabla_{Z_0} \Upsilon \rangle_q = \rho^2 \, L(w, \Phi) +
\rho \, Q(w, \Phi).
\]
Collecting these estimates we get
\[
\langle N, \nabla_{Z_0} \Upsilon \rangle_q  = \rho^2 \, L(w, \Phi)
+ Q(w, \Phi).
\]
Hence we conclude that
\[
\begin{array}{rlllll}
\langle N, \nabla_{Z_0} (\del_s w \, \Upsilon) \rangle_q & = &
\del_s^2 w \, \langle N, \Upsilon \rangle_q + \del_s w \, \langle
N, \nabla_{Z_0} \Upsilon \rangle_q \\[3mm]
&  = & - \del_s^2 \, w + Q(w,\Phi)
\end{array}
\]

\smallskip

\noindent {\bf Step 2} Next,
\[
\langle N, \nabla_{Z_0} \Phi' \rangle_q  = \rho \, \langle N ,
\Phi'' \rangle_q + \sum_{j=1}^n \del_{x_0} \phi_j \, \langle N,
\nabla_{Z_0} \, X_j \rangle_q
\]
where $\Phi''(x_0 ) :=  \sum \del_{x_0}^2 \phi_j (x_0)\,X_j$. From
(\ref{eq:3-9}), we have
\[
(\del_{x_0} \phi_j) \, \langle N, \nabla_{Z_0} X_j\rangle_q = \rho^2
\, L(w, \Phi) + \rho \, Q(w,\Phi).
\]
Also, using the same decomposition of $N$, and employing
(\ref{eq:3-8}) and (\ref{lems}),
\[
\begin{array}{rllll}
\langle N , \Phi'' \rangle_q & = & - \langle \Upsilon , \Phi''
\rangle_q + \langle N+\Upsilon, \Phi''\rangle_q \\[3mm]
& = & - \langle \Upsilon , \Phi'' \rangle_q + Q(w, \Phi)\\[3mm]
& = & - \langle \Upsilon , \Phi'' \rangle_p + \rho^2 \, L (w, \Phi) +
Q(w, \Phi)
\end{array}
\]
Collecting these gives
\[
\langle N, \nabla_{Z_0} \Phi' \rangle_q  = - \rho \, \langle \Upsilon ,
\Phi'' \rangle_p + \rho^2 \,L(w, \Phi) + \rho \, Q(w, \Phi).
\]

\smallskip

\noindent {\bf Step 3} Expanding $Z_0$ gives
\begin{equation}
\langle N, \nabla_{Z_0} X_0 \rangle_q = \rho \, \langle N,
\nabla_{X_0} X_0 \rangle_q +  \rho \, \del_s w \, \langle N,
\nabla_{\Upsilon} X_0 \rangle_q  + \rho \, \sum_{j=1}^n \phi_j \,
\langle N, \nabla_{X_j} X_0 \rangle_q \label{eq:3-11}
\end{equation}
With the help of (\ref{eq:3-7}) and (\ref{eq:3-8}), we evaluate
\[
\begin{array}{rllll}
\langle N, \nabla_\Upsilon X_0 \rangle_q & = & {\cal O} (\rho) + L (w, \Phi) + Q(w,\Phi)\\[3mm]
\langle N, \nabla_{X_j} X_0 \rangle_q & = & {\cal O} (\rho) + L (w, \Phi) + Q(w,\Phi)\\[3mm]
\langle N + \Upsilon , \nabla_{X_0} X_0 \rangle_q & = &  \rho \,  L (w,
\Phi) + Q(w,\Phi),
\end{array}
\]
and plugging these into (\ref{eq:3-11}) already gives
\[
\langle N, \nabla_{Z_0} X_0 \rangle_q = - \rho \, \langle \Upsilon ,
\nabla_{X_0} X_0 \rangle_q + \rho^2 \, L (w, \Phi) + \rho \,
Q(w,\Phi)
\]
Using (\ref{eq:3-2b}) in Lemma~\ref{le:3-1} we get the expansion
\[
\begin{array}{rllll}
\left. \nabla_{X_0} X_0\right|_q  & = & \ds \ds - \sum_{j=1}^n  \rho \,
\langle R(X_j, X_0) \, \Upsilon ,
X_0 \rangle_p \, X_j + {\cal O} (\rho^2)  \\[3mm]
& {} & -\, \ds \sum_{j=1}^n \langle R(X_j, X_0) \, \Phi , X_0 \rangle_p
\, X_j + \rho \, L(w, \Phi) + Q(w, \Phi),
\end{array}
\]
and so
\[
\begin{array}{rllllll}
\langle N, \nabla_{Z_0} X_0 \rangle_q & = & \ds \rho^2 \,
\sum_{j=1}^n \langle R(X_j, X_0) \Upsilon, X_0 \rangle_p \,
\langle \Upsilon ,  X_j \rangle_q  +  {\cal O} (\rho^3) \\[3mm]
& {} & +\,  \rho \, \sum_{j=1}^n \langle R(X_j, X_0) \,
\Phi, X_0 \rangle_p \, \, \langle \Upsilon , X_j \rangle_q \\[3mm]
& {} & +\,  \rho^2 \, L (w, \Phi) + \rho \, Q(w,\Phi).
\end{array}
\]
Finally, using (\ref{lems}) again, we conclude that
\[
\begin{array}{rlllll}
\langle N, \nabla_{Z_0} X_0 \rangle_q & = & \rho^2 \, \langle R(\Upsilon,
X_0 ) \,  \Upsilon , X_0 \rangle_p +
{\cal O} (\rho^3) +\rho \, \langle R(\Upsilon, X_0 ) \Phi , X_0 \rangle_p \\[3mm]
& {} & +\,  \rho^2 \, L(w, \Phi) + \rho \, Q(w, \Phi),
\end{array}
\]
which, together with the results of Step 1 and Step 2, completes
the proof of the first estimate.

\smallskip

\noindent
{\bf Estimate (\ref{eq:3-6b}):} Decompose
\[
\begin{array}{rllll}
\langle N, \nabla_{Z_0} Z_j \rangle_q & = & \rho \, \langle N, Y_j
\rangle_q \, \del_s w + \rho  \,  \langle N, \Upsilon\rangle_q \,
\del_s \del_{y_j} w \\[3mm]
& {} & +\, \rho \, (1+w) \, \langle N, \nabla_{Z_0} Y_j\rangle_q + \rho
\, \langle N, \nabla_{Z_0} \Upsilon\rangle_q.
\end{array}
\]
As above we use (\ref{eq:3-8}) and (\ref{eq:3-10}) to estimate
\[
\langle N, Y_j \rangle_q  = - \langle \Upsilon, Y_j\rangle_q  +
\langle N+ \Upsilon , Y_j\rangle_q = L(w, \Phi) + Q (w, \Phi)
\]
Similarly, by Lemma~\ref{le:3-2} and (\ref{eq:3-8}),
\[
\langle N, \Upsilon \rangle_q = - 1 + L(w, \Phi) + Q(w, \Phi)
\]
But now, by (\ref{eq:3-8}) and (\ref{eq:3-9}), we have
\[
\langle N, \nabla_{Z_0} \, Y_j \rangle_q = {\cal O}(\rho^2) + \rho
\, L(w, \Phi) + \rho  \, Q(w, \Phi)
\]
and
\[
\langle N, \nabla_{Z_0} \, \Upsilon \rangle_q = {\cal O}(\rho^2) + \rho
\, L(w, \Phi) + \rho  \, Q(w, \Phi),
\]
and the proof of the estimate follows directly.

\smallskip

\noindent {\bf Estimates (\ref{eq:3-6c}) and (\ref{eq:3-6d}):}
Observe that, thanks to the result of Proposition~\ref{pr:3-2}, we
can also write
\[
N = - \Upsilon + \frac{1}{\rho} \, \sum_{j=1}^n \alpha_j \, Z_j +
\hat N,
\]
where
\begin{equation}
\hat N = (L(w, \Phi) + Q(w, \Phi)) \, X_0 + \sum_{j=1}^n (\rho^2
\, L(w, \Phi) + Q(w, \Phi) ) \, X_j. \label{eq:hatn}
\end{equation}

Now write
\[
\begin{array}{rlllll}
\langle N, \nabla_{Z_j} Z_{j'} \rangle_q & = & \langle N,
\nabla_{Z_{j'}} Z_{j} \rangle_q \\[3mm]
& = & - \frac{1}{2} \, \left( \langle  \nabla_{Z_j}  N, Z_{j'}
\rangle_q + \langle  \nabla_{Z_{j'}}  N, Z_{j}
\rangle_q \right) \\[3mm]
& =& \frac{1}{2} \, \left( \langle  \nabla_{Z_j}  \Upsilon ,
Z_{j'} \rangle_q +  \langle  \nabla_{Z_{j'}}  \Upsilon ,
Z_{j} \rangle_q \right)\\[3mm]
& {} & -\, \frac{1}{2 \rho} \, \sum_{i=1}^{n} \,  \left( \langle
\nabla_{Z_j} (\alpha_i \, Z_i) , Z_{j'} \rangle_q + \langle
\nabla_{Z_{j'}} (\alpha_i \, Z_i) , Z_{j} \rangle_q \right) \\[3mm]
& {} & +\,  \frac{1}{2} \, \left( \langle \hat N, \nabla_{Z_j}
Z_{j'} \rangle_q + \langle \hat N, \nabla_{Z_{j'}} Z_{j}
\rangle_q \right) \\[3mm]
& {} & -\, \frac{1}{2} \, \left( \del_{y_j} \, \langle \hat N,
Z_{j'}\rangle |_q + \del_{y_{j'}} \, \langle \hat N, Z_{j}\rangle
|_q\right)
\end{array}
\]

\smallskip

\noindent {\bf Step 1} By (\ref{eq:3-7}), we can estimate
\[
\begin{array}{rlllll}
\nabla_{Z_j} Z_{j'} & = & \rho \, \del_{y_j} w \,
Y_{j'} + \rho \, \del_{y_j}\del_{y_{j'}} w \, \Upsilon \\[3mm]
& {} & +\, \rho \, (1+w) \, \nabla_{Z_j} Y_{j'} + \rho \,
\del_{y_{j'}}w \, \nabla_{Z_j} \Upsilon\\[3mm]
& = &  ({\cal O} (\rho^3) + \rho^2 \, L(w, \Phi) +
\rho^2 Q(w, \Phi) ) \, X_0  \\[3mm]
& {} & +\, \sum_{k=1}^n ({\cal O} (\rho^3) + \rho \, L(w, \Phi) +
\rho^2 \, Q(w, \Phi) ) \, X_k,
\end{array}
\]
Observe that the coefficient of $X_0$ is slightly better than the
coefficient of the other $X_k$ since the first two terms only
involve the $X_k$. Using this together with (\ref{eq:hatn}) we
conclude that
\[
\langle \hat N, \nabla_{Z_j} Z_{j'} \rangle_q + \langle \hat N,
\nabla_{Z_{j'}} Z_{j} \rangle_q =  \rho^3 \, L(w, \Phi) + \rho \,
Q(w, \Phi)
\]

\smallskip

\noindent {\bf Step 2} Next, using (\ref{eq:hatn}) together with
(\ref{lems}), we find that
\[
\del_{y_j} \, \langle \hat N, Z_{j'}\rangle |_q  + \del_{y_{j'}}
\, \langle \hat N, Z_{j}\rangle |_q = \rho^3 \, L(w, \Phi) + \rho
\, Q(w, \Phi)
\]

\smallskip

\noindent
{\bf Step 3} We now estimate
\[
A_{jj'} :=  \langle \nabla_{Z_j}  \Upsilon , Z_{j'} \rangle_q
+ \langle \nabla_{Z_{j'}}  \Upsilon , Z_{j} \rangle_q.
\]
It is convenient to define
\[
A'_{jj'}: = \frac{1}{1+w}\,\left(\langle
\nabla_{Z_j}(1+w)\,\Upsilon,Z_{j'}\rangle_q + \langle
\nabla_{Z_{j'}}(1+w)\,\Upsilon,Z_j\rangle_q \right) ,
\]
It follows from Lemma~\ref{le:3-2} that
\[
A_{jj'} = A_{jj'}' + + \rho \, Q(w,\Phi)
\] hence  it is enough to
focuss on the estimate of $A_{jj'}'$. To analyze this term, let us
revert for the moment and regard $w$ and $\Phi$ as functions of
the coordinates $(t,y)$ (rather than $(s,y)$), and also consider
$\rho$ as a variable instead of just a parameter. Thus we consider
\[
\tilde{F}(\rho,t,y) =  F \big(t, \rho(1+w(t,y))\Upsilon(y) +
\Phi(t)\big).
\]
The coordinate vector fields $Z_j$ are still equal to $\tilde F_*
(\del_{y_j})$, but now we also have $ (1+w)\Upsilon =  \tilde F_*
(\del_\rho)$, which is the identity we wish to use below. Now,
from (\ref{eq:3-4}), we write
\[
A'_{jj'} = \frac{1}{1+w} \, \left(\langle \nabla_{\del_\rho} Z_j ,
Z_{j'}\rangle_q + \langle
 \nabla_{\del_\rho}  Z_{j'} , Z_j \rangle_q\right) = \frac{1}{1+w} \,
 \del_\rho \langle Z_j, Z_{j'} \rangle |_q
\]
Therefore, it follows from (\ref{eq:3-4}) in
Proposition~\ref{pr:3-2} that
\[
\begin{array}{rlllll}
A_{jj'} & = &  \frac{1}{1+w} \, \del_\rho \, [ \rho^2 \, \langle
Y_j, Y_{j'}\rangle_p
+ \frac{\rho^4}{3} \, \langle R(\Upsilon, Y_j) \, \Upsilon, Y_{j'}
\rangle_p + {\cal O} (\rho^5)\\[3mm]
& {} & +\, 2 \, \rho^2 \,w \, \langle Y_j, Y_{j'}\rangle_p + \ds
\frac{\rho^3}{3}   \, (\langle R(\Upsilon, Y_j)\, \Phi , Y_{j'}
\rangle_p + \langle
R(\Upsilon, Y_{j'}) \, \Phi , Y_j \rangle_p) \\[3mm]
& {} & +\, \rho^4 \, L(w, \Phi) + \rho^2 \, Q(w, \Phi) ] + \rho \, Q(w, \Phi) \\[3mm]
& = & \frac{1}{1+w} \,[ 2 \, \rho \,  \langle Y_j, Y_{j'}\rangle_p
+ \frac{4}{3} \, \rho^3 \, \langle R(\Upsilon, Y_j) \, \Upsilon,
Y_{j'} \rangle_p + {\cal O} (\rho^4)\\[3mm]
& {} & \, 4 \, \rho \,w  \, \langle Y_j, Y_{j'}\rangle_p  + \rho^2 \,
\left( \langle R(\Upsilon, Y_j)\, \Phi, Y_{j'} \rangle_p  + \langle
R(\Upsilon, Y_{j'})\, \Phi, Y_j \rangle_p \right) \\[3mm]
& {} & +\, \rho^3 \, L(w, \Phi) + \rho \, Q(w, \Phi) ] \\[3mm]
& = & \ds 2 \, \rho \, \langle Y_j, Y_{j'}\rangle_p  +
\frac{4}{3} \, \rho^3 \, \langle R(\Upsilon, Y_j) \, \Upsilon, Y_{j'} \rangle_p
+ {\cal O} (\rho^4)\\[3mm]
& {} & +\, \ds  2 \, \rho \,w  \, \langle Y_j, Y_{j'}\rangle_p + \rho^2
\, (\langle R(\Upsilon, Y_j) \, \Phi , Y_{j'} \rangle_p +
\langle R(\Upsilon, Y_{j'})\, \Phi , Y_{j} \rangle_p) \\[3mm]
& {} & +\, \rho^3 \, L(w, \Phi) + \rho \, Q(w, \Phi)
\end{array}
\]

\noindent {\bf Step 4} Finally, we must compute
\[
\begin{array}{rllll}
B_{jj'} & : = & \langle  \nabla_{Z_j} (\alpha_i \, Z_i) , Z_{j'}
\rangle_q  + \langle  \nabla_{Z_{j'}} (\alpha_i \, Z_i) ,
Z_{j} \rangle_q  \\[3mm]
& =  & \langle Z_i , Z_{j'} \rangle_q \, \del_{y_j} \alpha_i +
\langle Z_i , Z_j \rangle_q \, \del_{y_{j'}} \alpha_i \\[3mm]
& {} & +\,\alpha_i \, ( \langle  \nabla_{Z_j} Z_i , Z_j' \rangle_q
 + \langle  \nabla_{Z_{j'}} Z_i , Z_j \rangle_q )\\[3mm]
& =  & \langle Z_i , Z_{j'} \rangle_q \, \del_{y_j} \alpha_i +
\langle Z_i , Z_j \rangle_q \, \del_{y_{j'}} \alpha_i \\[3mm]
& {} & +\, \alpha_i \, ( \langle  \nabla_{Z_i} Z_j , Z_j' \rangle_q
 + \langle  \nabla_{Z_i} Z_{j'} , Z_j \rangle_q )\\[3mm]
& =  & \langle Z_i , Z_{j'} \rangle_q \, \del_{y_j} \alpha_i +
\langle Z_i , Z_j \rangle_q \, \del_{y_{j'}} \alpha_i + \alpha_i
\, \del_{y_i} \,  \langle  Z_j , Z_{j'} \rangle_q
\end{array}
\]
Observe that, (\ref{eq:3-55}) implies
\[
\del_{y_j} \langle Y_i, Y_{j'} \rangle_p  = 0.
\]
Using this together with (\ref{eq:3-4}) and the expression for the
$\alpha_i$ given in Proposition~\ref{pr:3-3}, we get
\[
\alpha_i \, \del_{y_i} \,  \langle  Z_j , Z_j' \rangle_q  = \rho^4
\, L(w, \Phi) + \rho^2 \, Q(w, \Phi)
\]

It follows from (\ref{eq:3-4}) and the definition of $\alpha_i$
again that
\[
\langle Z_i , Z_{j'} \rangle_q \, \del_{y_j} \alpha_i = \rho^2 \,
\langle Y_i , Y_{j'} \rangle_p \, \del_{y_j} \alpha_i  + \rho^4 \,
L(w, \Phi) + \rho^2 \, Q(w, \Phi)
\]
Therefore, it remains to estimate $\langle Y_i , Y_{j'} \rangle_p
\, \del_{y_j} \alpha_i$. By definition, we have
\[
\sum_{i=1}^n \alpha_i \, \langle Y_i, Y_{j'}\rangle_p =
 \del_{y_{j'}} w + \frac{\rho}{3} \, \langle R(\Phi, \Upsilon)
\, \Upsilon, Y_{j'}\rangle_p
\]
Differentiating with respect to $y_j$ we get
\begin{equation}
\sum_{i=1}^n \left( \langle Y_i, Y_{j'} \rangle_p \, \del_{y_j}
\alpha_i + \alpha_i \, \del_{y_j} \langle Y_i, Y_{j'} \rangle_p
\right) =  \del_{y_j} \del_{y_{j'}} w + \frac{\rho}{3} \,
\del_{y_j} \langle R(\Phi, \Upsilon) \, \Upsilon, Y_{j'}\rangle_p
\label{eq:start}
\end{equation}
Again, it follows from (\ref{eq:3-55}) that $\del_{y_j} \langle
Y_i, Y_{j'} \rangle_p  = 0$. Moreover, using (\ref{eq:3-9}), we
first estimate
\[
\nabla_{Z_j} \Upsilon =  Y_j + {\cal O}(\rho^2) + \rho \, L(w,
\Phi) +  \rho \, Q(w,\Phi);
\]
and, using in addition (\ref{eq:3-55}), we also get
\[
\nabla_{Z_{j}} Y_{j'} = a \,\Upsilon + {\cal O}(\rho^2) + \rho \,
L(w, \Phi) + \rho \, Q(w,\Phi)
\]
for some $a \in \R$. Reinserting this in (\ref{eq:start}) yields
\[
\begin{array}{rlllll}
\sum_{i=1}^n \langle Y_i, Y_{j'} \rangle_p  \, \del_{y_j} \alpha_i
& = & \del_{y_j } \del_{y_{j'}} w + \frac{\rho}{3} \,
\langle R(\Phi, Y_j) \, \Upsilon, Y_{j'} \rangle_p  \\[3mm]
& {} & +\, \frac{\rho}{3} \, \langle R(\Phi, \Upsilon) \, Y_j , Y_{j'}
\rangle_p + \rho^3 \, L(w,\Phi) + \rho^2 \, Q(w,\Phi),
\end{array}
\]
since $R(\Upsilon, \Upsilon) =0$.

Collecting these estimates, we conclude that
\[
B_{jj} = \rho^2 \, \del_{y_j}^2 w \, + \frac{\rho^3}{3} \, \langle
R(\Phi, Y_j ) \, \Upsilon, Y_j \rangle_p  + \rho^4 \, L(w, \Phi) +
\rho^2 \, Q(w, \Phi)
\]
since $ \langle R(\Phi, \Upsilon) \, Y_j , Y_j \rangle_p = 0$ and
also that
\[
B_{jj'} = \rho^2 \, L(w, \Phi) + \rho^2 \, Q(w, \Phi)
\]
when $j \neq j'$. With the estimates of the previous steps, this
finishes the proof of the last two estimates! \hfill $\Box$

\subsection{The mean curvature}
Collecting the estimates of the last subsection and taking the
trace, we have now proved that the mean curvature $H(w, \Phi)$ of
the hypersurface $\calT_\rho(w , \Phi)$ satisfies
\begin{eqnarray}
& \ds n \, \rho \, H(w, \Phi) - (n-1) \   =    \notag \\
& \left( \frac{2}{3} \, \langle R ( \Upsilon , \, X_0) \, \Upsilon ,
X_0 \rangle_p  - \frac{1}{3} \,\mbox{Ric} (\Upsilon, \Upsilon) \, \right) \,
\rho^2 + {\cal O} (\rho^3) \notag \\
& \ds  - \left( \del_s^2 w + \Delta_{S^{n-1}} w + (n-1) \, w \right)
-  \rho \, \langle \Phi''+  R(\Phi, X_0) \, X_0 ,\Upsilon\rangle_p
\label{eq:mc} \\
& \ds +  \rho^2 \, L(w, \Phi) + Q(w, \Phi). \notag
\end{eqnarray}
(We recall that if $E_\alpha$ is an orthonormal basis of $T_pM$,
then
\[
\mbox{Ric}(\Upsilon,\Upsilon) := -\sum_{\alpha=0}^n \langle
R(\Upsilon,E_\alpha)\Upsilon, E_\alpha\rangle_p.)
\]

Denote by $(\lambda_j,\varphi_j)$ the eigendata of
$\Delta_{S^{n-1}}$, where the eigenfunctions are orthonormal and
counted with multiplicity.

A most important observation is that the second and third terms in the
expansion of $n \, \rho \, H$ are quadratic in the coordinates $x_j$.
Hence, when $\Phi=w=0$, we have
\begin{equation}
\begin{array}{rcl}
\left(n\,\rho\, H-(n-1),\varphi_j\right)_{L^2(S^{n-1})} & = & \calO(\rho^3),
\qquad j = 1, \ldots, n, \\[3mm]
\left(n \, \rho\,H - (n-1),\varphi_j\right)_{L^2(S^{n-1})} &= & \calO(\rho^2),
\qquad j \neq 1, \ldots, n,
\end{array}
\label{eq:exvan}
\end{equation}
or in other words, writing $f = n\rho H - (n-1)$, then $f =
\calO(\rho^2)$ but its $L^2 (S^{n-1})$ projection over $\varphi_1,
\ldots, \varphi_n$ satisfies $\hat{f} = \calO(\rho^3)$.

\section{Jacobi operators}

In this section we examine the mapping properties of some linear
operators which appear in the expression of the mean curvature of
$\calT_\rho (w, \Phi)$ given in (\ref{eq:mc}).

\subsection{Definitions}
The two linear operators appearing in the third line of (\ref{eq:mc}) are
\begin{equation}
\begin{array}{rlllll}
w & \longmapsto &  \calL_{SN\Gamma} \, w& := & \del_s^2 \, w +
\Delta_\theta \, w + (n-1) \, w, \\[3mm]
 \Phi & \longmapsto &\frakJ \, \Phi & := & \nabla_{X_0}^2\Phi +
R(\Phi,X_0)\,X_0.
\end{array}
\end{equation}
The latter is the Jacobi operator on $\Gamma$ corresponding to the
second variation of the length functional on curves, while (up to
a multiplicative factor) the former is the Jacobi operator for the
second variation of the area functional about a Euclidean cylinder
${\R}\times S^{n-1}(\rho)$.

Recall that the geodesic $\Gamma$ is said to be nondegenerate
when $\frakJ$ is invertible, i.e.\ if the equation $\frakJ \,
\Phi = 0$ has no nontrivial solutions on all of $\Gamma$. For a
generic metric on $M$, it is well known that all closed geodesics
are nondegenerate.

On the other hand, since it is already naturally expressed in terms of the
scaled coordinate $s$,
\[
\calL_{SN\Gamma} :\calC^{2, \alpha}_\rho (SN\Gamma) \longrightarrow
\calC^{0,\alpha}_\rho (SN\Gamma)
\]
is bounded uniformly in $\rho$. We can analyze this operator using
the eigendecomposition for $\Delta_\theta$ on $S^{n-1}$. As in \S
2.2, if the eigenfunction decomposition of $w$ is given by
\[
w (s, \theta) = \sum_{j\geq 0} \, w_j (s) \, \varphi_j (\theta),
\]
then $w$ decomposes as $w_0 + \hat{w} + \tilde{w}$, where
\[
\hat w : = \sum_{j=1}^n w_j \, \varphi_j \qquad \mbox{and} \qquad
\tilde w : = \sum_{j>n} w_j \, \varphi_j.
\]
We denote by $\Pi_0$, $\hat{\Pi}$ and $\tilde{\Pi}$ the
projections on to these three components, respectively. From now
on, we assume that we are working with functions $w$ such that
$\hat{\Pi} \, w = 0$, and thus we only need to be concerned with
the operators $(\calL_{SN\Gamma})_0$ and $\tilde{\mathcal
L}_{SN\Gamma}$ induced on the two other components. Note in
particular that
\[
(\calL_{SN\Gamma})_0 := \del_s^2 + n-1.
\]

\subsection{Mapping properties}

We now study the mapping properties of $\frakJ$ and (the
components of) $\calL_{SN\Gamma}$.

We first note that
\[
\frakJ : \calC^{2,\alpha}(\Gamma,N\Gamma)  \longrightarrow
\calC^{0,\alpha}(\Gamma,N\Gamma)
 \]
is an isomorphism when $\Gamma$ is a nondegenerate geodesic.

Next, we also assert that
\[
\tilde{\calL}_{SN\Gamma}: \tilde{\Pi} \,
\calC^{2,\alpha}_\rho(SN\Gamma) \longrightarrow \tilde{\Pi} \,
\calC^{0,\alpha}_\rho(SN\Gamma)
\]
is an isomorphism with inverse uniformly bounded as $\rho \to 0$;
this follows from the fact that $\Delta_\theta + (n-1) \leq -C <
0$ on this subspace. Details are left to the reader.

Finally, it is clear that
\[
(\calL_{SN\Gamma})_0: \Pi_0 \, \calC^{2,\alpha}_\rho(SN\Gamma)
\longrightarrow \Pi_0 \, \calC^{0,\alpha}_\rho(SN\Gamma)
\]
is bounded for every $\rho>0$, but is only invertible when
\[
\sqrt{n-1}\,\frac{\Lambda}{\rho} \notin 2\, \pi \, \Z;
\]
in the exceptional cases, there is a two-dimensional nullspace
spanned by
\[
\cos(\sqrt{n-1}\,s), \quad \sin(\sqrt{n-1}\,s),
\]
and hence
a two-dimensional cokernel. To determine the norm of its inverse
when $\sqrt{n-1}\,\Lambda/\rho \notin 2\,\pi\,\Z$, suppose that
$(\calL_{SN\Gamma})_0 \, v = f$. Then
\[
\begin{array}{rllll}
\sqrt{n-1} \, v(s) & = & \ds \sin (\sqrt{n-1} \, s ) \, \left(
\alpha + \int_0^{s} \cos (\sqrt{n-1} \, \sigma) \, f(\sigma)
\, d\sigma \right) \\[3mm]
& - & \ds \cos(\sqrt{n-1} \, s)\, \left( \beta +
\int_0^s\sin(\sqrt{n-1} \, \sigma)\,f(\sigma)\,d\sigma \right)
\end{array}
\]
where the constants $\alpha,\beta$ are chosen so that $w$ is
$\Lambda /\rho$-periodic. We find
\[
|\alpha| + |\beta|\leq \frac{c}{\rho
\,(1-\cos(\sqrt{n-1}\Lambda/\rho)) } \, ||f||_{L^\infty
(SN\Gamma)}
\]
for some constant $c >0$, independent of $\rho$, and from this we have
\begin{equation}
||v||_{\calC^{2,\alpha}_\rho (SN\Gamma)}\leq c \, \left(
||f||_{\calC^{0,\alpha}_\rho (SN\Gamma)} +
\frac{1}{\rho\,(1-\cos(\sqrt{n-1}\Lambda/\rho))}
\,||f||_{L^\infty(SN\Gamma)}\right), \label{eq:estimee1}
\end{equation}
again for some constant $c>0$ independent of $\rho$. Note that
when $f\in\calC^1$, there is an equivalent formula
\[
\begin{array}{rllll}
(n-1) \, v(s) - f (s) & = & \sin (\sqrt{n-1} \, s) \, \left(
\alpha - \int_0^s\sin(\sqrt{n-1}\, \sigma)\,
\del_s f(\sigma)\,d\sigma \right) \\[3mm]
& + & \cos (\sqrt{n-1}\,s)\, \left( \beta -
\int_0^s\cos(\sqrt{n-1}\, \sigma)\,\del_s f(\sigma) \, d\sigma
\right),
\end{array}
\]
where again $\alpha, \beta$ are chosen so that $w$ is
$\Lambda/\rho$\, - periodic. We now obtain
\[
|\alpha| + |\beta|\leq \frac{c}{\rho
\,(1-\cos(\sqrt{n-1}\Lambda/\rho)) } \, ||\del_s f||_{L^\infty
(SN\Gamma)}
\]
for some constant $c >0$ independent of $\rho$, so that
\begin{equation}
||v||_{{\cal C}^{2, \alpha}_\rho (SN\Gamma)} \leq  c \, \left(
\,||f||_{\calC^{0, \alpha}_\rho (SN\Gamma)} + \frac{1}{\rho \, (1-
\cos (\sqrt{n-1} \Lambda / \rho))} \, ||\del_s f||_{L^\infty
(SN\Gamma)}) \right). \label{eq:estimee2}
\end{equation}

\section{The constant mean curvature foliation}
We now use the results of \S 2 and \S 3 to perturb
$\calT_\rho(\Gamma)$ to a constant mean curvature hypersurface, at
least for $\rho$ sufficiently far from values where $({\cal
L}_{SN\Gamma})_0$ is degenerate.

According to the analysis of \S 2, we must find $w \in \calC^{2,
\alpha}_\rho (SN\Gamma)$ and $\Phi \in \calC^{2,
\alpha}(\Gamma,N\Gamma)$ such that
\begin{equation}
n\,\rho \, H(w, \Phi) = n-1
\label{eq:cmceq}
\end{equation}

Let us denote by
\[
f := \ds \frac{2}{3} \, \langle R ( \Upsilon , \, X_0) \Upsilon , X_0
\rangle_p  \, \rho^2 - \frac{1}{3} \, \mbox{Ric} (\Upsilon,
\Upsilon) \, \rho^2 + {\cal O} (\rho^3),
\]
the inhomogeneous term appearing in (\ref{eq:mc}) which corresponds to the
mean curvature when $w=\Phi = 0$. As usual, this decomposes into three
components, $f_0 + \hat{f} + \tilde{f}$, where $\hat{f}$ corresponds to a
section of the normal bundle which we write as $\rho \Psi$.
We are searching for $w=w_0 + \tilde{w}$ and $\Phi$ which satisfy the
coupled system
\begin{equation}
\left\{
\begin{array}{clllll} (\calL_{SN\Gamma})_0 \, w_0  & = & f_0 +
\rho^2 \, L(w,\Phi) + Q(w,\Phi) \\[3mm]
\frakJ \,  \Phi  & = &  \Psi + \rho \, L(w,\Phi) + \frac{1}{\rho} \, Q(w, \Phi) \\[3mm]
\tilde{\calL}_{SN\Gamma} \tilde{w}  & = & \tilde{f} + \rho^2 \,L(w, \Phi) + Q(w,\Phi) \\
\end{array}
\right.
\label{eq:system}
\end{equation}

We use the function space
\[
{\mathcal E}^{2,\alpha}_\rho : = \Pi_0 \,
\calC^{2,\alpha}_\rho(SN\Gamma) \oplus
\calC^{2,\alpha}(\Gamma,N\Gamma) \oplus \tilde{\Pi} \,
\calC^{2,\alpha}_\rho(SN\Gamma),
\]
where, for $\Xi = (w_0,\Phi,\tilde{w})$,
\[
\|\Xi \|_{{\mathcal E}^{2,\alpha}_\rho} : =  (1 -
\cos(\sqrt{n-1}\Lambda/\rho)) \,
\|w_0\|_{\calC^{2,\alpha}_\rho(SN\Gamma)}+
\|\Phi\|_{\calC^{2,\alpha}(\Gamma,N\Gamma)} +
\|\tilde{w}\|_{\calC^{2,\alpha}_\rho(SN\Gamma)}.
\]

The linear operators appearing on the left in (\ref{eq:system}) are all
invertible provided $\sqrt{n-1}\,\Lambda /\rho \notin \Z$. Thus,
multiplying by their inverses, we rewrite this system as
\[
\Xi = {\mathfrak N}(\Xi),
\]
and so we solve our problem by finding a fixed point of ${\mathfrak N}$
in ${\mathcal E}^{2,\alpha}_\rho$.

\begin{lemma} Write $\Xi_0 : = {\mathfrak N}(0)$. Then
for $\sqrt{n-1}\,\Lambda/\rho\notin \Z$, we have
\[
\|\Xi_0\|_{{\mathcal E}^{2,\alpha}_\rho} \leq \frac{c_0}{2}\, \rho^2
\]
for some $c_0 > 0$.
\end{lemma}
{\bf Proof~:} Clearly
\[
\|f_0\|_{{\calC}^{0,\alpha}_\rho(SN\Gamma)}+ \rho^{-1} \, ||\del_s
f_0||_{L^\infty(SN\Gamma)} \leq c\, \rho^2;
\]
moreover
\[
||\tilde{f}||_{{\calC}^{0,\alpha}_\rho(SN\Gamma)} \leq c\,\rho^2
\qquad \mbox{and} \qquad ||\Psi||_{{\calC}^{0,\alpha}(\Gamma,N\Gamma)}\leq c \, \rho^2.
\]
Since, by definition,
\[
\Xi_0 = ((\calL_{SN\Gamma})_0^{-1} \, f_0, \frakJ^{-1}\,\Psi_0,
(\tilde{\calL}_{SN\Gamma})^{-1} \,\tilde{f}),
\]
the result follows from (\ref{eq:estimee2}) and the uniform
boundedness of the inverses of the inverses of these linear
operators as $\rho \to 0$. \hfill $\Box$

\smallskip

Next, from the properties of the operators $L$ and $Q$ we deduce the
\begin{lemma}
There exists a constant $c >0$ such that, for the same $c_0$ as in
the previous Lemma, and for any $\Xi_1, \Xi_2 \in {\mathcal
E}^{2,\alpha}_\rho$ satisfying
\[
\|\Xi_i\|_{{\mathcal E}^{2, \alpha}_\rho} \leq c_0 \, \rho^2,
\]
we have
\[
\|{\mathfrak N}(\Xi_2) - {\mathfrak N}(\Xi_1)\|_{{\mathcal
E}^{2,\alpha}_\rho} \leq c  \,
\frac{\rho}{(1-\cos(\sqrt{n-1}\Lambda/\rho))^2} \, \|\Xi_2
-\Xi_1\|_{{\mathcal E}^{2, \alpha}_\rho}.
\]
\end{lemma}
{\bf Proof~:} It follows from (\ref{eq:AAA}) that
\[
\| \rho^2 \, L (w, \Phi)\|_{{\mathcal C}^{0,\alpha}_\rho} \leq c
\, \frac{\rho^2}{1-\cos(\sqrt{n-1}\Lambda/\rho)} \,
\|\Xi\|_{{\mathcal E}^{2, \alpha}_\rho}
\]
if $w=  w_0+ \tilde w$ and $\Xi = (w_0, \Phi , \tilde w)$.
Moreover, if $\|\Xi_i\|_{{\mathcal E}^{2, \alpha}_\rho} \leq c_0
\, \rho^2$, then we have from (\ref{eq:BBB})
\[
\| Q (w_2, \Phi_2) - Q (w_1, \Phi_1)\|_{{\mathcal
C}^{0,\alpha}_\rho} \leq c \,
\frac{\rho^2}{(1-\cos(\sqrt{n-1}\Lambda/\rho))^2} \, \|\Xi_2
-\Xi_1\|_{{\mathcal E}^{2, \alpha}_\rho}.
\]
where $w_i = w_{0,i}+ \tilde w_i$ and $\Xi = (w_{0,i}, \Phi_i ,
\tilde w_i)$. Now, the result follows at once from the inequality
\[
\|\Phi\|_{{\calC}^{0, \alpha} (\Gamma, N\Gamma)} \leq \, c \,
\rho^{-\alpha} \, \| \Phi\|_{{\calC}^{0, \alpha}_\rho (\Gamma,N\Gamma)}
\]
and the uniform bounds on $\tilde{\calL}_{SN\Gamma}^{-1}$ 
and ${\mathfrak J}^{-1}$, and the bound on
$({\cal L}_{SN\Gamma})_0^{-1}$ given in (\ref{eq:estimee1}).
Details are left to the reader. \hfill $\Box$

\smallskip

Collecting these results, we now have the
\begin{proposition}
Fix  $\alpha  \in (0,1)$. Then there exists a $c_1>0$ such that if
$k$ is sufficiently large and $\rho$ satisfies
\[
\frac{1}{k+1} + \frac{c_1}{k^{9/4}} \leq \frac{2\pi}{\sqrt{n-1}
\,\Lambda} \,\rho \leq \frac{1}{k} - \frac{c_1}{k^{9/4}}
\]
then there exists a solution $(w_0, \tilde w ,\Phi)$ of
(\ref{eq:system}) in ${\mathcal E}^{2,\alpha}_\rho$. This solution
satisfies
\[
\|(w_0, \tilde w ,\Phi)\|_{{\mathcal E}^{2,\alpha}_\rho} \leq c_0
\, \rho^2
\]
\label{pr:tre}
\end{proposition}
{\bf Proof~:} It is easy to check that
\[
\frac{\rho^{1-\alpha} }{(1-\cos(\sqrt{n-1}\Lambda/\rho))^2}
\]
is as small as we want, provided $c_1$ is chosen large enough. It
is then easy to check that, when $k$ is large enough, ${\mathfrak
N}$ is a contraction from the ball of radius $c_0 \, \rho^2$ into
itself. \hfill $\Box$

\smallskip

This proposition yields the existence of CMC perturbations of the
tubes $\calT_\rho(\Gamma)$ for all radii $\rho \in I_k$, when $k$
is large. We shall denote the perturbation functions as $w_\rho$
and $\Phi_\rho$ to emphasize their dependence on $\rho$. We shall
now revert to thinking of these as depending on $x_0$ rather than $s$,
and in particular we write
\[
w_\rho (x_0,\theta) :=  w_0 (x_0/\rho) + \tilde w(x_0/\rho,\theta).
\]
Following through the proof, it is not hard to see that these functions
depend smoothly on $\rho$. Furthermore, since the tubes $\calT_\rho(\Gamma)$
already foliate, it suffices to verify that the mapping
\begin{equation}
(\rho,x_0,\theta) \longmapsto G(x_0, \rho(1+w_\rho)\theta + \Phi_\rho)
\label{eq:diffeo}
\end{equation}
is a local diffeomorphism.

First,
\[
||w_\rho||_{L^\infty(SN\Gamma)} + ||\Phi||_{L^\infty(\Gamma,N\Gamma)}
\leq c \,\rho^2.
\]
Also, from the construction itself, we have
\[
||\del_\rho w_\rho||_{L^\infty (SN\Gamma)} + ||\del_\rho
\Phi||_{L^\infty (\Gamma, N\Gamma)} \leq c \, \rho.
\]
These certainly imply that (\ref{eq:diffeo}) is a local\
diffeomorphism, and hence our CMC surfaces form a local foliation.
\hfill $\Box$

\section{Explaining the gaps}

In special cases, such as when $\Gamma$ is a circle in the flat
torus $T^{n+1} = S^1 \times \R^{n}/ a \Z^n$, all of the geodesic
tubes about $\Gamma$, of any radius, have constant mean curvature,
and thus there are no gaps in the local foliation. On the other
hand, Theorem~1 only provides for a local foliation with gaps. As
indicated in the introduction, there are good reasons why this
construction fails from working at all radii. We explain this in
greater detail now. We first show that for generic metrics, the
moduli space of CMC surfaces isotopic to a geodesic tubes
$\calT_\rho(\Gamma)$ in $M \setminus \Gamma$ is smooth and
one-dimensional. The index of the Jacobi operator is constant
along components, and by estimating the index for the surfaces
close to $\calT_\rho(\Gamma)$, $\rho \in I_k$, we show that there
are infinitely many components of this moduli space. We conclude
by examining in more detail a very degenerate case, where all of
the geodesic tubes $\calT_\rho(\Gamma)$ are CMC (when smooth), so
there are no gaps, but we prove that in this situation there are
infinitely many bifurcating branches.

\subsection{The moduli space}

Denote by $\calH(M,\Gamma,g)$ the moduli space of all CMC surfaces
$\Sigma \hookrightarrow M$ which are isotopic to any one of the
geodesic tubes $\calT_\rho(\Gamma)$, for $\rho$ small, in
$M\setminus\Gamma$, with respect to the metric $g$.

\begin{proposition}
There is an open dense set $\calU$ of metrics (in the
$\calC^{m,\alpha}$ topology for any $m \geq 3$) on $M$ such that
for $g \in \calU$, $\calH(M,\Gamma,g)$ is a smooth one-dimensional
manifold.
\end{proposition}
Fix a surface $\Sigma_0$ in the correct isotopy class, which has
CMC with respect to some metric $g_0$. Nearby surfaces may be
written as  normal graphs over $\Sigma_0$, and hence are
parametrized by (small) scalar functions on $\Sigma_0$. Now
consider the mapping
\[
G: \calC^{m,\alpha}(M,S^2T^*M) \times \calC^{m,\alpha}(\Sigma_0) \longrightarrow
\calC^{m-2,\alpha}(\Sigma_0)
\]
which assigns to a metric $g$ and a scalar function $w$ on
$\Sigma_0$ the mean curvature function of the submanifold
$\Sigma_w = \{x + w(x)\nu(x) \, : \,  x \in \Sigma_0\}$, regarded
as a function on $\Sigma_0$, with respect to $g$. It is not hard
to show \cite{Whi} (and also \cite{Maz-surv}) that the
differential of this mapping is always surjective, and moreover
the restriction of this differential to the tangent space of the
second factor (which is simply the Jacobi operator) is Fredholm of
index zero, and is an isomorphism except when there exist
nontrivial Jacobi fields. The result is then a straightforward
application of the Sard-Smale theorem, since $G$ is transverse to
the one-dimensional curve of constant functions in the range
space.

On the other hand, applying the construction of Theorem 1 when the
metric $g \in \calU$, we obtain a set of smooth  one-dimensional
families CMC surfaces in $\calH(M,\Gamma,g)$. Let $G_g(w) =
G(g,w)$; an implication of this proposition is that when $G_g$ is
regular at $w$ and the mean curvature of $\Sigma_w$ is equal to
$H$, then some interval $(H-\epsilon, H+\epsilon)$ parameterizes
$\calH(M,\Gamma,g)$ locally near $\Sigma_w$, or in other words,
the mean curvatures of the CMC surfaces near to $\Sigma_w$ assume
all values near to $H$. Although in this case the corresponding
surfaces form a local foliation, this need not be true in general.

\subsection{The index}
We now claim that for generic $g$, $\calH(M,\Gamma,g)$ has
infinitely many components. In the following, for $\rho \in I_k$,
let $\Sigma_\rho$ denote the CMC hypersurface constructed
in Theorem~1.
\begin{proposition}
Let $g$ be a metric for which $\calH(M,\Gamma,g)$ is a smooth
one-dimensional manifold. Then for each sufficiently large value
of $k$, the surfaces
$\Sigma_\rho$, $\rho \in I_k$, lie in different components of
$\calH(M,\Gamma,g)$ and have index equal to $\mbox{Index} \, (\Gamma ) + 2 \, k + 1$.
\end{proposition}

We prove this theorem by computing the index of $\Sigma_\rho$;
by definition this is the number of negative eigenvalues of the Jacobi
operator ${\mathcal L}_\rho$ on $\Sigma_\rho$. This index is locally
constant in $\calH(M,\Gamma,g)$ when this moduli space is nondegenerate,
but since we shall show that the index increases with $k$, this
will imply the result.

It follows from (\ref{eq:mc}) and the properties of solution
$(w_0, \tilde w, \Phi)$ given in Proposition~\ref{pr:tre} that the
Jacobi operator about $\Sigma_\rho$, i.e.\ the linearization
of the operator $(w_0,\tilde{w},\Phi) \mapsto n\rho H$, has the form
\begin{equation}
\begin{array}{rcl}
\ds \calL_\rho ( v,\Psi)  & = &  -\left( \del_s^2 v +
\Delta_{S^{n-1}} \, v + (n-1) \, v \right) \ds -\rho \, \langle
\Psi''+ R(X_0,
\Psi) X_0, \theta \rangle_p \\[3mm] & {}& \ + \ \rho^2 L(v, \Psi),
\end{array}
\label{eq:jacop}
\end{equation}
where the (linear operator) $L$ satisfies the usual assumptions.

\begin{lemma}
The quadratic form associated to $\calL_\rho$ has the expansion
\[
\begin{array}{rlllll}
B_\rho (v,\Psi) & : = & \int_{SN\Gamma} (|\del_s v|^2 + |\nabla_\theta v|^2
- (n-1) \, v^2 ) \\[3mm]
& {} & \, + \, \omega_n \, \rho \, \int_{\Gamma} (|\del_t \Psi|^2 - \langle R (\Psi, X_0) \,
\Psi,  X_0 \rangle_p )\\[3mm]
& {} & \ + \, \rho^2 \, C (v, \Psi)
\end{array}
\]
where $v \in \Pi_0 H^1_\rho (SN\Gamma) \oplus \tilde \Pi \, H^1_\rho(SN\Gamma)$
(i.e.\ the Sobolev space $H^1(SN\Gamma)$ defined with respect to the vector
fields $\del_s$ and $\del_{y_j}$), $\Psi \in H^1 (\Gamma, N\Gamma)$,
$\omega_n$ is a positive constant depending only on the dimension,
and where $C$ is a quadratic form satisfying
\[
|C(v, \Psi)|\leq c \, \left( \int_{SN\Gamma} (|\del_s v|^2 + |\nabla_\theta
v|^2 + v^2 ) + \int_{\Gamma} (|\del_t \Psi|^2 + |\Psi |^2) \right)
\]
for some constant $c >0$ independent of $\rho$.
\end{lemma}
{\bf Proof~:} For $v$ and $\Psi = \sum \psi_j E_j$ as in this statement, the
fibrewise linear function on $SN\Gamma$ corresponding to $\Psi$ is $\hat{v}
= \sum \psi_j \theta_j$. We can either regard $B_\rho$ as a quadratic form
in the variables $(v,\Psi)$ or in $v + \hat{v}$. Then, by definition
\[
B_\rho(v+\hat{v}) = \int_{SN\Gamma}(\calL_\rho(v,\Psi))(v+\hat{v}).
\]

Inserting the expression (\ref{eq:jacop}) and integrating over $\Sigma_\rho$,
we obtain the first summand in the expression for $B_\rho$, involving only $v$,
without difficulty. Next, integrating over the spherical fibres of $SN\Gamma$, we have
\[
-\int_{SN\Gamma} \langle \Psi'', \theta\rangle \langle \Psi, \theta\rangle =
\sum_{j,k=1}^n \int_{SN\Gamma} \psi_j' \psi_k' \theta_j \theta_k = \omega_n
\sum_{j=1}^n \int_{\Gamma} (\psi_j')^2;
\]
there is a similar reduction for the term in $\frakJ$ of order $0$ to an
integral over $\Gamma$.
The error term leading to $C$ arises from the error term in
$\calL_\rho$, as well as the discrepancy in this last calculation
caused by using the volume form on $SN\Gamma$ rather than the one on
$\Sigma_\rho$. The first of these error terms,
\[
\int L (v, \Psi) \, (v+ \hat{v}),
\]
is almost of the correct form. However, since $\calL_\rho$ involves
$\Psi''$, this error might include terms of the form $v \, \Psi''$,
which are at first glance too big since, integrating by parts, they
equal $\rho^{-1} \, \del_s v \, \Psi'$. However, examining the computations
leading up to (\ref{eq:mc}), one can check that $\Phi''$ only enters
through the term
\[
\rho \, \langle N,\Phi''\rangle_p (1+\calO(\rho^2)+ L(w,\Phi) + Q(w,\Phi))
\]
(recall that we compute $\rho$ times the mean curvature), which has
linearization
\[
\rho \, \langle N, \Psi''\rangle_p (1+\calO(\rho^2)+ L(w,\Phi)+
Q(w,\Phi)) =\rho\langle\Psi'',\theta\rangle_p + \calO(\rho^3)\, L (\Psi).
\]
Hence this gives, at worst, terms like $\rho^2 \, \del_sv \, \Psi'$.
It is much more straightforward to check that all the other terms in $C$
satisfy the correct bounds. \hfill $\Box$

\smallskip

As usual, write $v(s, y) = v_0(s) + \tilde v(s,y)$, where both summands
are orthogonal to the linear eigenfunctions $\varphi_j$, $1 \leq j \leq n$;
we also identify $\hat{v}(t,y) :=  \langle \Psi,\theta\rangle_p$.
Thus for $v = v_0 + \hat v + \tilde v \in \Pi_0 H^1_\rho (SN\Gamma)
\oplus \hat\Pi \, H^1 (SN\Gamma) \oplus \tilde \Pi \, H^1_\rho
(SN\Gamma)$ we have the quadratic form
\[
B_\rho(v) = B_0 (v_0)+\rho \, \hat B(\hat v) + \tilde B(\tilde v) +
\rho^2 \, C(v),
\]
where
\[
\begin{array}{rlllll}
B_0(v_0) & : = & \omega_n \,  \int_\Gamma (|\del_s v_0|^2 - (n-1) \, v^2_0 ) \\[3mm]
\hat B(\hat v) & : = &  \omega_n \, \int_\Gamma (|\del_t \Psi|^2 -
\langle R (\Psi,
X_0) \, \Psi, X_0\rangle_p )\\[3mm]
\tilde B( \tilde v) & : = &  \int_{SN\Gamma} (|\del_s \tilde v|^2 +
|\nabla_\theta \tilde v|^2 - (n-1) \, \tilde v^2 )
\end{array}
\]

Assuming that $\rho \in I_k = (\rho'_k , \rho''_k)$, let us now
compute the index of $B$. Since this index is locally constant in
$\rho$, we shall choose
\[
\rho = \frac{\sqrt{n-1}}{4\, \pi} \, \Lambda \,
\left(\frac{1}{k}+ \frac{1}{k+1} \right),
\]
which is directly in the middle of $I_k$.

Clearly,
\[
\int |\del_s \tilde v|^2 +|\nabla \tilde v|^2  + \tilde v^2 \leq c
\, \tilde B (\tilde v).
\]

Next, decompose $\hat v = \hat v^+ + \hat v^-$, where $\hat v^\pm$
lies in the sum of the eigenspaces of $\frakJ$ with positive or negative
eigenvalues, respectively. (Recall that this operator is assumed to be
nondegenerate, hence has no zero eigenspace.) We then have
\[
\int |\del_t \hat v^+|^2 + |\hat v^+|^2 \leq  c \, \hat B(\hat
v^+)
\]
and
\[
\int |\del_t \hat v^-|^2 + | \hat v^- |^2 \leq  - c \, \hat B(\hat v^-).
\]

We can similarly decompose the remaining component $v_0$ as
$v_0^+ + v_0^-$, where $v_0^\pm$ lie in the eigenspaces corresponding
to the positive or negative eigenvalues of $(\calL_{SN\Gamma})_0 =
-\del_s^2 - (n-1)$. Using that $\rho$ lies in the middle of $I_k$,
we can estimate
\[
\int |\del_s v_0^+|^2 + |v_0^+|^2 \leq  c \, \rho^{-1} \, B_0(v_0^+)
\]
and
\[
\int |\del_s v_0^-|^2 + |v_0^- |^2 \leq  - c \, \rho^{-1} \,  B_0(v_0^-).
\]

Using all of these estimates, we now obtain that
\[
\begin{array}{rllll}
| C (v) | & \leq & c\, \rho^{-1} \, \left( B_0 (v_0^+) - B_0(v_0^-) \right)
\\[3mm]
& {} & +\,  c \, \left( \hat B (\hat v^+) - \hat B (\hat v^-) \right)  +
c \, \tilde B (\tilde v)
\end{array}
\]
for some constant $c>0$ independent of $\rho$. This gives, in turn,
\[
B' (v) \leq B_\rho(v)  \leq B'' (v)
\]
where
\[
\begin{array}{rlllll}
B''(v) & = &  (1+ c\, \rho) \, \left( B_0 (v_0^+) +
\rho \, \hat B (\hat v^+) \right) \\[3mm]
& {} & +\, (1- \, c \, \rho) \, \left( B (v_0^-) + \rho \,  \hat B (\hat v^-)
 \right)\\[3mm]
 & {} & +\, (1+c \, \rho^2) \, \tilde B (\tilde v)
\end{array}
\]
and
\[
\begin{array}{rlllll}
B'(v) & = &  (1- c\, \rho) \, \left( B_0 (v_0^+) +
\rho \, \hat B (\hat v^+) \right) \\[3mm]
& {} & +\, (1 + \, c \, \rho) \, \left( B (v_0^-) + \rho \, \hat B (\hat v^-) \right)
 \\[3mm]
& {}& +\, (1-c \, \rho^2) \, \tilde B (\tilde v).
\end{array}
\]

These upper and lower bounds on $B$ imply that
\[
\mbox{Index}  \, (B'') \leq \mbox{Index} \, (B_\rho) \leq
\mbox{Index} \, ( B').
\]
The proof is completed by the following
\begin{lemma}
When $\rho$ is small enough, the index of $B'$ and $B''$ are both equal
to $\mbox{Index} \, (\Gamma) + 2\, k+1$.
\end{lemma}
{\bf Proof~:} If $\rho$ is chosen so that $1- c \, \rho > 1/2$, then
the index of $B'$ and $B''$ equals the sum of the dimensions of the
spaces on which $B_0$ and $\hat B$ are negative. But these equal
$2k+1$ and $\mbox{Index} \, (\Gamma)$, respectively. \hfill $\Box$

\subsection{Bifurcations in a degenerate case}
We consider in more detail the (very) degenerate case where
$(M,g)$ is the flat torus $T^{n+1} = S^1 \times \R^{n}/a \Z^n$ for
$a$ sufficiently large, and $\Gamma = S^1 \times \{0\}$. After
moding out by all the continuous symmetries, the moduli space is
still one-dimensional, but has infinitely many singularities.

Each of the geodesic tubes $\calT_\rho(\Gamma) := \Sigma_\rho =
S^1 \times S^{n-1} (\rho)$ is CMC, with mean curvature $H_\rho =
(n-1)/\rho$. We only consider the case where $\Sigma_\rho$ is
embedded, i.e.\ when $\rho < 1/2$. The Jacobi operator for
$\Sigma_\rho$ is
\[
L_\rho = \Delta_{\Sigma_\rho} + |A_{\Sigma_\rho}|^2.
\]
In terms of our standard cylindrical coordinates, $t \in S^1$, $\theta \in S^{n-1}$,
\[
L_\rho = \del_t^2 + \frac{1}{\rho^{2}}\left(\Delta_\theta + (n-1)\right).
\]
Introducing eigendata $\{\phi_k(t),k^2\}$ and $(\psi_\ell(\theta),
-\lambda_\ell^2)$ in each component, we see that $L_\rho$ reduces
to multiplication by $ B(k,\ell,\rho) = -k^2 +
\rho^{-2}(n-1-\lambda_\ell^2)$ on the $(k,\ell)$ eigenspace. Since
$\lambda_\ell^2$ is always of the form $j(n-2+j)$ for some
nonnegative integer $j$, and hence $n-1-\lambda_\ell^2 \leq 0$
unless $\ell=0$. This gives the
\begin{proposition} The surface $\Sigma_\rho$ is always degenerate;
its nullspace consists of the span of the eigenmodes
$\phi_0(t)\psi_\ell(\theta)$, $\ell = 1, \ldots, n$ and, in case
$\rho^2 = (n-1)/k^2$ for some $k \in \N$, also
$\phi_k(t)\psi_0(\theta)$.
\end{proposition}
The `trivial degeneracies' are those comprised by the first set of
elements, which exist for all $\rho$. These correspond to the
obvious geometric fact that translating $\Sigma_\rho$ parallel to
itself in any direction normal to $\Gamma$ gives a family of CMC
surfaces with the same mean curvature. These can  be eliminated if
we mod out by these symmetries. Namely, using linear coordinates
$x = (x_0,x_1, \ldots, x_n)$ in $T^{n+1}$, let $G$ be the finite
group generated by reflections in the $x_j = 0$ plane, $j = 1,
\ldots, n$.
\begin{corollary} Acting on the space of $G$-invariant functions on $SN\Gamma$,
the Jacobi operator $L_\rho$ is degenerate if and only if $\rho^2
= (n-1)/k^2$, $k \in \N$.
\end{corollary}
These degeneracies have a direct geometric explanation too, for it
is precisely at these radii where the family of `$k$-bump'
Delaunay surfaces begins to develop.

We now have the picture that $\calH(T^{n+1},S^1,g_0)$ consists of
the union of the interval $(0,a/2)_\rho$ and infinitely many other
intervals $(0,(n-1)/k^2]_{\e_k}$, where the variable $\e_k$ is the
Delaunay necksize in the $k^{\mathrm{th}}$ bifurcating branch. In
other words, the moduli space looks like an open interval with
infinitely many spines sticking out of it. When the metric on
$T^{n+1}$ is perturbed generically, this moduli space smooths out;
the singularity at each degenerate radius disappears, and this
`spiny interval' breaks into infinitely many components. The CMC
surfaces for these slightly perturbed metrics are small
perturbations either of the geodesic tubes or else of the Delaunay
surfaces, except near the turning points. The gaps encountered in
our construction correspond exactly to the small regions around
these degenerate radii where the moduli space is curving away from
the interval $(0,a/2)_\rho$.

\section{Limits of constant mean curvature foliations}
There is a sort of converse to Theorem~\ref{th:existence} which we
can prove regardless of the dimension of $\Gamma$. Let
$\Gamma$ be a closed $\ell$-dimensional submanifold of $M^{n + 1}$,
$1 \leq \ell \leq n$, and suppose that there exists a sequence of
hypersurfaces $\Sigma_j$ such that
\begin{itemize}
\item[a)] $\Sigma_j$ has constant mean curvature $(n-\ell)/n\rho_j$, where $\{\rho_j\}$
is a decreasing sequence with $\rho_j \searrow 0$;
\item[b)] $\Sigma_j$ is isotopic in $M \setminus \Gamma$ to the tube
$\calT_{\rho_j}(\Gamma)$;
\item[c)] There exists a $c > 0$ such that $\Sigma_j$ is contained
inside $\calT_{c\rho_j}(\Gamma)$ for all $j$.
\end{itemize}
Item c) implies that $\Sigma_j \to \Gamma$ in Hausdorff distance; note
also that we are not requiring the existence of any sort of local foliation,
just a sequence of CMC hypersurfaces converging to $\Gamma$. We conjecture
that some set of hypotheses very near to these (for example, assuming also
a bound on the second fundamental form, as in the introduction, or that $\Sigma_j$
is trapped between the tubes of radius $c'\rho_j$ and $c\rho_j$ for fixed constants
$0 < c' < c$.) should be enough to ensure that $\Gamma$ is minimal.
Unfortunately, this seems to be quite difficult to prove, and so we shall
restrict ourselves to a very special situation by making a fourth, quite
restrictive, hypothesis:
\begin{itemize}
\item[d)] There exist functions $w_j \in \calC^{2,\alpha}(SN\Gamma)$,
$\Phi_j \in \calC^{2,\alpha}(\Gamma,N\Gamma)$ satisfying
\[
||w_j||_{\calC^{2,\alpha}} + ||\Phi_j||_{\calC^{2,\alpha}} \leq c \rho_j^2
\]
for some $c>0$, independent of $\rho_j$, such that
\[
\Sigma_j = \calT_{\rho_j}(w_j,\Phi_j).
\]
\end{itemize}
Quite important (and potentially restrictive) here is that the norms of
$w_j$ and $\Phi_j$ are bounded in $\calC^{2,\alpha}$, {\it not} $\calC^{2,\alpha}_\rho$.

\begin{theorem}
Let $\Gamma$ be a $\calC^2$ compact embedded $\ell$-dimensional
submanifold of $M$ for some $1 \leq \ell \leq n$,  and suppose
that $\Sigma_j$ is a sequence of CMC hypersurfaces converging to
$\Gamma$ and satisfying the hypotheses a) -- d). Then $\Gamma$ is
minimal.
\end{theorem}

The proof is based on an argument from geometric measure theory
which is now fairly standard in the analysis of such `condensation
problems', cf.\ \cite{Bet-Orl-Sme-1}, \cite{Bet-Orl-Sme-2}.

We drop the subscript $j$ and consider a functional for which each
$\Sigma_\rho$ is critical. The argument proceeds by writing the formula which
expresses the fact that the first variation of this functional vanishes,
and then taking the limit of this formula (in a very weak sense) as $\rho \to 0$.
The limiting first variation equation implies the minimality of $\Gamma$.

Any one of the CMC hypersurfaces $\Sigma = \Sigma_\rho$ bound a
compact domain $D(\Sigma)$, which is the component of $M\setminus
\Sigma$ containing $\Gamma$. Define the measure
\[
d\mu_\Sigma = dA_\Sigma - n \, H \, dV_\Sigma,
\]
where $H$ is the (constant) mean curvature of $\Sigma$ and where we have set
\[
dA_{\Sigma} : = \calH^n\lfloor_{\Sigma} \qquad \mbox{and}
\qquad  dV_\Sigma : = \calH^{n+1} \lfloor_{D(\Sigma)}
\]
($\calH^k$ is $k$ dimensional Hausdorff measure).

CMC hypersurfaces are critical for the functional $\Sigma \to \int d\mu_\Sigma$.
In other words, if $X$ is any $\calC^2$ vector field on $M$, and $\phi_t$
the associated one-parameter family of diffeomorphisms, then
\begin{equation}
\left. \int d\mu_{\phi^*_t \Sigma} \right|_{t=0} = 0.
\label{eq:crit}
\end{equation}
We now compute this variation another way.  In fact, for any hypersurface $\Sigma$
and any continuous function $f$, we derive that
\[
\left. \del_t \int f \, dA_{\phi^*_t \Sigma} \right|_{t=0} =
\int X f \, \, dA_\Sigma + \int f\,(\mbox{div}X-\langle\nabla_N X,N\rangle)\,dA_\Sigma,
\]
where $N$ is the unit normal to $\Sigma$, and similarly,
\[
\left. \del_t \int f \, dV_{\phi^*_t\Sigma} \right|_{t=0} =
\int X f \, dV_\Sigma + \int f \, \mbox{div} X \, dV_\Sigma.
\]
Hence, setting $f \equiv 1$ and using (\ref{eq:crit}), we obtain
\begin{equation}
0 = \left. \del_t \int f d\mu_{\phi^*_t \Sigma}\right|_{t=0} =
\int \, \mbox{div} X \, d\mu_\Sigma - \int \langle \nabla_N X, N\rangle \, dA_\Sigma.
\label{eq:crit2}
\end{equation}

Next, let $dL_\Gamma = \calH^\ell\lfloor_{\Gamma}$ and denote by $\omega_{n-\ell}$
the volume of $S^{n-\ell}$ with respect to its standard metric (thus
$\omega_{n-\ell}/(n+1-\ell)$ is the volume of the ball $B^{n+1-\ell}$).
\begin{lemma}
As $\rho \searrow 0$,
\[
\rho^{\ell-n}\,dA_{\Sigma_\rho}\rightharpoonup \omega_{n-\ell} \,dL_\Gamma,
\qquad \mbox{and}\qquad
\rho^{\ell-n-1}\, dV_{\Sigma_\rho} \rightharpoonup \frac{\omega_{n-\ell}}{n+1-\ell}\,dL_\Gamma
\]
in the sense of measures. In particular,
\[
\rho^{\ell-n}\,d\mu_{\Sigma_\rho}\rightharpoonup
\frac{\omega_{n-\ell}}{n+1-\ell} \, dL_\Gamma
\]
\label{le:limits1}
\end{lemma}
{\bf Proof~:} This follows from Fubini's theorem and the fact that the functions
$w$ and $\Phi$ appearing in the parametrization of $\Sigma$ are uniformly
controlled in $\calC^{2,\alpha}$ as $\rho \to 0$.
\hfill $\Box$

\smallskip

On the other hand, we also have
\begin{lemma}
Let $E_1, \ldots, E_{n+1-\ell}$ be a local orthonormal frame for $N\Gamma$.
Then as $\rho \searrow 0$,
\[
\rho^{\ell-n} \, \langle \nabla_N X, N\rangle_p \,
dA_{\Sigma_\rho} \rightharpoonup \frac{\omega_{n-\ell}}{n+1-\ell}
\, \sum_i \langle \nabla_{E_i} X,  E_i \rangle_p \, dL_\Gamma
\]
in the sense of measures.
\label{le:limits2}
\end{lemma}
{\bf Proof~:} As before, using Fubini's theorem and the uniform control
on $w$ and $\Phi$, it suffices to check that this formula holds for
the round sphere of radius $\rho$ in $\R^{n+1}$ as $\rho \to 0$, where
again the formula is a standard computation.
\hfill $\Box$

\smallskip

Now multiply (\ref{eq:crit2}) by $\rho^{\ell-n}$ and let $\rho \to 0$. From these
two lemmas, we conclude that
\[
\int  \left( \mbox{div} X - \sum_i \langle \nabla_{E_i}  X , E_i
\rangle \right) \, dL_\Gamma = 0
\]
On the other hand, if $F_1, \ldots, F_\ell$ is a local orthonormal frame for $T\Gamma$,
then
\[
\mbox{div} X = \sum_{j=1}^{\ell} \langle \nabla_{F_j}X, F_j\rangle + \sum_{i=1}^{n+1-\ell}
\langle \nabla_{E_i}X, E_i\rangle,
\]
so this last equation is equivalent to
\[
\int  \sum_j  \langle X , \nabla_{F_j} F_j \rangle \, dL_\Gamma = 0.
\]
Since the vector field $X$ is arbitrary, we conclude that the normal
component of  $\sum_j \nabla_{F_j} F_j$ is equal to $0$. This implies that
the mean curvature of $\Gamma$ vanishes, i.e.\  that $\Gamma$ is minimal.
\hfill $\Box$

\smallskip

As already discussed at the beginning of this section, it would be
much nicer to prove this theorem under less stringent hypotheses.
We conclude by discussing this in more detail.

Suppose that there exist sequences of intervals $I_k =
(\rho_k',\rho_k'')$ in $\R^+$ with with $\rho_k',\rho_k'' \to 0$,
such that for each $k$ and $\rho \in I_k$ there exists a CMC
hypersurface $\Sigma_\rho$ isotopic to $\calT_\rho(\Gamma)$ in $M
\setminus \Gamma$. Suppose furthermore that this hypersurface
satisfies:
\begin{itemize}
\item[a')] the hypersurfaces $\{\Sigma_\rho\}_{\rho \in I_k}$ form a local foliation;
\item[b')] the mean curvature of $\Sigma_\rho$ equals $\frac{n-\ell}{n}\,\frac{1}{\rho}$;
\item[c')] there exists a constant $c >0$, independent of $k$ and $\rho \in I_k$,
such that
\[
\Sigma_\rho \subset B_\rho(\Gamma) := \{q \in M^{n+1}: \quad \mbox{dist}_g
(q,\Gamma)\leq c\,\rho\}.
\]
\item[d')] There exists a constant $c >0$, again independent of $k$ and $\rho \in I_k$,
such that $|A_{\Sigma_\rho} |\leq c/\rho$.
\end{itemize}
We conjecture that these hypotheses alone are sufficient to
conclude that $\Gamma$ is minimal. Indeed, it is possible to prove
many of the necessary facts, but a few crucial ones seem much more
difficult to obtain.

If we rescale $\Sigma_\rho$ from a point $p \in \Gamma$ by the
factor $1/\rho$, we obtain a family of CMC surfaces which are
cylindrically bounded. We can obtain area bounds for these
rescaled hypersurfaces, just as in \cite{Ye-1}, and so conclude
that at least along subsequences, this family converges to a
complete embedded cylindrically bounded hypersurface in
$\R^{n+1}$. It is known \cite{Kus-Kor-Mee} (and also
\cite{Maz-surv}) that all such hypersurfaces must lie in the
family of Delaunay unduloids $D_\e$; however, using that the
$\Sigma_\rho$ are leaves of a local CMC foliation, we obtain a
global bounded positive Jacobi field on the limiting surface, and
it may be checked directly that in the Delaunay family, only the
cylinder admits such a Jacobi field. Thus far we have shown that,
having fixed $p \in \Gamma$ and rescaling about this point, then
some subsequence of these surfaces, say $\Sigma_j$, converges to a
cylinder with axis parallel to the rescaled limit of $\Gamma$. The
first difficulty is a standard one in the subject: if we knew the
uniqueness of this limit, then we could straightaway conclude the
existence of functions $w_j$ and $\Phi_j$ such that $\Sigma_j =
\calT_{\rho_j}(w_j,\Phi_j)$. We would also be able to conclude
that
\[
\|w_j \|_{\calC^{2,\alpha}_\rho} + \|\rho^{-1}_j \, \Phi_j \|_{\calC^{2,\alpha}_\rho} = o (1).
\]
However, we only obtain these bounds in $\calC^{2,\alpha}_\rho$, not
$\calC^{2,\alpha}$. Because of this, we are unable to obtain bounds on
the volume of $\Sigma$ of the form
\[
\int_{\Sigma_{\rho}} dA_\Sigma \leq C \rho^{n-\ell};
\]
on the other hand, it is quite straightforward to prove that
$\int_{D(\Sigma_\rho)} dV_{\Sigma} \leq C \rho^{n+1-\ell}$. In any
event, we are only able to show that the second conclusion of
Lemma~\ref{le:limits1} holds, and so we are unable to take the
weak limit of $d\mu_{\Sigma_\rho}$.  The final difficulty arises
in proving the analogue of Lemma~\ref{le:limits2}, and this is the
case because of the rather weak control we have for the
derivatives of $w$ and $\Phi$ in the $\Gamma$ direction.

Plausibly, the most general theorem of this sort would involve a
sequence of CMC hypersurfaces which are known to converge in Hausdorff
distance, and perhaps satisfying hypotheses a') -- d'). The conclusion
should be that $\Sigma_\rho$ converge to a minimal submanifold
$\Gamma$, of some dimension $\ell$, away from a set of Hausdorff measure
smaller than $\ell$. For example, if $\ell=1$, it is quite conceivable
that such a set of surfaces might converge to a broken geodesic (satisfying
certain constraints at the break points).

\end{document}